%% file: AASpectralGapTheoremSU_d1.tex
\input amstex
\documentstyle{amsppt}
\magnification 1200
\UseAMSsymbols
\hsize 5.5 true in
\vsize 8.5 true in
\parskip=\medskipamount
\NoBlackBoxes
\def\span{\text { span\,}}
\def\Tr{\text{Tr\,}}
\def\mathbb{\Bbb}

\def\mathcal{\Cal}

\def\supp{\text{\rm supp\,}}

\def\Mat{\text{\rm Mat}}
\def\Res{\text{\rm Res}}
\def\Proj{\text{\rm Proj}}

\def\ve{\varepsilon}
\def\vp{\varphi}
\def\arrowk{^\to{\kern -6pt\topsmash k}}
\def\arrowK{^{^\to}{\kern -9pt\topsmash K}}
\def\arrowr{^\to{\kern-6pt\topsmash r}}
\def\bark{\bar{\kern-0pt\topsmash k}}
\def\arrowvp{^\to{\kern -8pt\topsmash\vp}}
\def\arrowf{^{^\to}{\kern -8pt f}}
\def\arrowg{^{^\to}{\kern -8pt g}}
\def\arrowu{^{^\to}a{\kern-8pt u}}
\def\arrowt{^{^\to}{\kern -6pt t}}
\def\arrowe{^{^\to}{\kern -6pt e}}
\def\tk{\tilde{\kern 1 pt\topsmash k}}
\def\barm{\bar{\kern-.2pt\bar m}}
\def\barN{\bar{\kern-1pt\bar N}}
\def\barA{\, \bar{\kern-3pt \bar A}}

\def\mathbb{\Bbb}

\def\nint{\not\!\!\int}
\def\dist{\text{\rm dist\,}}

\def\diam{\text {\rm \,diam\,}}
\TagsOnRight
\NoRunningHeads
\document
\lineskip=1pt\baselineskip=12pt\lineskiplimit=0pt
\def\Im{\text {Im\,}}

\def\supp{\text{\rm supp\,}}

\def\Mat{\text{\rm Mat}}
\def\Re{\text{\rm Re}}

\def\ve{\varepsilon}
\def\vp{\varphi}
\hsize = 6.2true in

\def\diam{\text {\rm \,diam\,}}
\vsize=8.2 true in
\TagsOnRight
\NoRunningHeads
\topmatter
\title
A Spectral Gap Theorem in $SU(d)$
\endtitle

\author
J.~Bourgain, A.~Gamburd
\endauthor
\address
School of Mathematics, Institute for Advanced Study, Princeton, NJ 08540
\endaddress
\email bourgain\@ias.edu
\endemail
\address
Mathematics Ph.D. Program,  The CUNY Graduate Center, 365 Fifth Avenue, New York, NY 10016
\endaddress
\email agamburd\@gc.cuny.edu
\endemail
\thanks
{The first author was supported, in part,  by NSF grants DMS-0808042 and DMS-0835373.  The second author was supported, in part, by DARPA, via AFOSR grant
FA9550-08-1-0315, by NSF grant DMS-0645807, and by the Sloan Foundation}
\endthanks
\endtopmatter

{\bf Abstract:}  We establish the spectral gap property for dense subgroups of  $SU(d)$  ($d\geq 2$), generated by  finitely many elements with algebraic entries; this result was announced in [BG3].
The method of  proof differs,  in several crucial aspects,  from that used in [BG]  in the case of $SU(2)$.

\vskip 15 pt

\centerline
{\bf \S0. Introduction and Outline}

For $k \ge 2$ let  $g_1, \dots, g_k$ be a finite set of elements
in $G=SU(d)$  ($d\geq 2$).   We associate with them an averaging (or Hecke) operator
$z_{g_1, \dots, g_k}$,  taking $L^2(SU(d))$
into $L^2(SU(d))$: $$ z_{g_1, \dots, g_k} f(x) =
\sum_{j=1}^{k}(f(g_j x) +f(g_j^{-1}(x)).$$  We denote by
$\supp(z)$ the set $\{g_1, \dots, g_k, g_1^{-1}, \dots,
g_k^{-1}\}$ and by $\Gamma_{z}$ the group generated by $\supp(z)$.
It is clear that $z_{g_1, \dots, g_k}$ is self-adjoint and that
the constant function is an eigenfunction of $z$ with eigenvalue
$\lambda_0(z)=2k$. Let $\lambda_{1}(z_{g_1, \dots, g_k})$ denote
the supremum of the eigenvalues of $z$ on the orthogonal
complement of the constant functions in $L^2(SU(d))$. We say that
$z$
 has  a spectral gap  if
$$\lambda_{1}(z_{g_1, \dots, g_k})<2k.$$

It is common to, alternatively, refer to the situation described above, by asserting that the spectral gap property holds for  $\Gamma_{z}$.

In this paper we generalize the result  on the spectral gap for finitely generated subgroups  of  $SU(2)$, established in [BG],  to dense subgroups of $SU(d)$ ($d \geq 2$),  generated by finitely many elements with algebraic entries.

\proclaim
{Theorem 1}  Assume that  $\{g_1, \ldots, g_k\}\in SU(d) \cap \Mat_{d\times d} (\overline{\Bbb Q})$,  and  that the group generated by
$g_1, \ldots, g_k$ is  Zariski dense\footnote{Note that the Zariski density assumption is equivalent to the topological density of the group generated by  $\{g_1, \dots, g_k\}$ in  $SU(d)$.} in $SL_d(\Bbb C)$.
Then the associated  Hecke  operator $z_{g_1, \dots, g_k}$  has a spectral gap.
\endproclaim

Various applications of such a spectral gap result   (to, among other things, the Banach-Ruziewicz problem,  the theory of  quasi-crystals, and the question, arising in the theory of quantum computation,  of whether a ``computationally universal''  set is necessarily ``efficiently universal'')  are discussed
in [BG].  

It should be pointed out, however, that the method of proof in the present paper differs, in several crucial aspects,  from the one given  in [BG]  in the case of $SU(2)$.

In [BG], the proof of the spectral gap proceeded by, first, establishing a ``product theorem''  for general subsets of $SU(2)$.  Both the statement and the proof of the latter result is  not  unrelated
to the product theorem in  $SL_2(p)$, established by Helfgott [H]  (and generalized to the groups of higher rank by Breuillard, Green and Tao  [BGT] and  by Pyber and Szabo  [PS]\footnote{Mention should be made, too, of the groundbreaking work of Hrushovski [Hr], and of the work of Breuillard and Green [BrGr] on the classification of approximate subgroups of the unitary group, yielding elegant and far-reaching generalization of Jordan's theorem [J].});
a key ingredient in the proof of the pertinent  product theorem  in the aforementioned papers   is the exact size of intersections of multiplicatively stable subsets of the group with maximal tori.

In contrast, the approach we follow  in the present paper is akin to the one in [BG1, 2], 
and  is 
based, crucially, on multi-scale arguments (available for groups defined over  $\Bbb C$ or $\Bbb Z/p^n\Bbb Z$),  and Lie algebra point of view\footnote{It might be worth remarking, that,  in the stressed crucial reliance on multi-scale and Lie algebra structures, this approach is  reminiscent of the Solovay-Kitaev algorithm in quantum computation [DN].}.
The salient features of this approach can be encapsulated as follows:  (a) first, using ``tools from arithmetic combinatorics'', we construct in the ``approximate group''  (see [Tao],  [BG], [BG1, 2] for background) 
``approximately one-dimensional structure''  in a suitable neighborhood of the identity;  (b) subsequently, this structure 
 serves  as the main  building block to recover the full Lie algebra; (c)
certain ``escape'' (from hyperplanes)  issues,  coming into play in (b),   are addressed,  using, in an essential way,  the theory of random matrix products.

In connection with (c), it should be pointed out,  that,  whereas in [BG2]  the ``classical''  theory of random matrix products  (see, for example,  [BL])   for Zariski-dense
subgroups of $SL_d$ (as developed by, among others,  Furstenberg [F], Goldsheid-Margulis [GM],  and Guivarch [G])  was directly applicable, 
in the present $SU(d)$ setting nontrivial difficulties   arise,  due to the absence of proximal (in the obvious sense)  elements,  necessitating the use of  non-Archimedian local fields and exterior powers of the Lie algebra (cf. [A]).

It is our expectation, that  the method of proof of the spectral gap result developed in the present paper in the context of $SU(d)$, should also be  applicable to other continuous semi-simple  Lie groups; we  intend to pursue this in a forthcoming  paper.

The main ingredient from arithmetic combinatorics, alluded to in (a) above,  is an extension of the ``discretized ring theorem'' (see [B1, 2] ) from $\Bbb R $ to $\Bbb C$ (see Proposition 2 at the end of the Introduction), and, crucially, to  Cartesian products $\Bbb C^d$.
This extension is obtained in \S1--\S8: parts of the argument are closely related to  [B2]; 
several steps are presented  in a somewhat greater generality than what is, strictly speaking, necessary for the purposes of this paper.

Returning to the proof of the spectral gap result in $SU(d)$, let us conclude this introduction by giving a rough summary of  the various steps in the argument.

As mentioned above, the overall approach is akin  to the one used in [BG1, 2].
Let
$$
\nu= \frac 1{2k} \sum^k_{i=1} (\delta_{g_i}+\delta_{g_{i}^{-1}})\tag 0.1
$$
be the probability measure supported on the generators $g_1, \dots, g_k$.
   Denoting by $P_\delta$ an approximate identity on $G=SU(d)$, and taking $\delta\to 0$, our first objective is to prove that, for $\tau>0$ (some fixed small constant), we can
ensure that
$$
\Vert\nu^{(\ell)}* P_\delta\Vert_\infty <\delta^{-\tau}\tag 0.2
$$
for $\ell< C(\tau) \log \frac 1\delta$.
Here $\nu^{(\ell)}$ refers to the $\ell$-fold convolution power of $\nu$.
This is achieved by iterating an ``$L^2$-flattening lemma'' (see lemma 10.7 in \S10) which is the main technical step in this part of the argument.
First, an application  (originating in  [BG0], and, by now,  standard) of  non-commutative Balog-Szemeredi-Gowers  lemma (proved in [Tao]),  reduces the matter  to the study of  ``approximate groups''  $H\subset G$.
Note that these objects are defined combinatorially, and, a priori, have no algebraic structure.
Our goal is to show, roughly speaking, that if $H$ is a $\delta$-approximate group such that $\nu^{(\ell)}(H)$ is ``large''   (where $\ell \sim \log \frac 1\delta$), then $H$ has to
be ``almost all''   of $G$, up to $\delta$-approximation.
This will, then,  provide the desired contradiction.

The first step in our program  is to produce in $H$ a large set of elements that are  ``approximately diagonal''
(in a suitable basis).  The key idea,  underlying the proof in this step, originates in the  work of Helfgott [H].
Let us point out, however, that, in contrast  to [H],  and to the  subsequent work on the product theorems  in $SL_d(p)$ and other finite simple groups 
([BGT], [PS]), the precise size of our diagonal set is not important.
The construction of this almost diagonal set appears in \S9.
The fact that the generators $g_1, \ldots, g_k$  have  algebraic entries plays a role here, but not the assumption on the Zariski density  of  $\langle g_1, \dots, g_k \rangle$. 

The relevant statement is Proposition 9 in \S9.  It should be stressed, that,  compared with [BG1, 2]  (and, also,  with   [BGT], [PS]), a significant difference is that we do not rely on regular elements to produce the almost diagonal set, and Proposition 9 provides such  construction in a greater generality.

The ``almost diagonal'' set of matrices is processed further using the discretized ring theorem in $\Bbb C^d$, resulting 
 in our main building block:  a structured,  ``essentially one-dimensional''  set in the Lie algebra.
A further amplification requires addressing  certain ``escape'' issues that depend on the assumption of Zariski density.
Thus, in
\S11, we establish the crucial ``$L^2$-flattening  lemma''
for convolution powers, conditional on the ``escape from hyperplane'' assumption (*), which is addressed in \S12 for $d=3$
and \S14 in the general case.

Similarly to the approach originating in  [BG1, 2],    proving the escape property relies on the theory of random matrix products.
Recall that the two main assumptions in this theory are  proximality   and  strong irreducibility.
In contrast to the case of Zariski dense subgroups of $SL_d(\Bbb R)$, elements in $SU(d)$ are obviously not expanding in the usual sense,  and the application of 
the  theory of random matrix products  requires considering non-Archimedean places (here, again, we use the fact that the elements $g_1, \ldots, g_k$ are algebraic) and, also,
representation  on wedge-products of the adjoint representation (see \S 12, 13, 14).
A treatment of random matrix product theory in the context of general local fields may be found in the recent paper of Aoun  [A].

Once (0.2) is established, the final step in the proof of the spectral gap, requires application of basic results pertaining to the representation theory of $SU(d)$. 
One way to proceed  (as was done in [BG]),  is to use the idea, originating in the work of   Sarnak and Xue [SX], of exploiting ``high multiplicity'' of
nontrivial egenvalues (which follows from ``high dimensionality''
of nontrivial irreducible representations);  in the continuous
setting of a compact group this idea was  implemented in [GJS]  by summing  over
the suitably chosen range of representations, and then applying Poisson
summation. 
In \S10, we follow a different route,   which,  in a sense, is more ``geometric'' (cf.  [BY]).
First, a new argument for $d=2$ is given.
Next, using $SU(2)$-subgroups in $SU(d)$, the general case is treated (this type of argument was used earlier in the work of Burger and Sarnak [BS]).
One of the ingredients is a convolution principle, stated in Lemma 10.35, which  appears to be a rather basic result, of independent interest,  pertaining to the  harmonic analysis on  the unitary group.

The paper is subdivided into two parts.  In the second part  (\S9--\S14), Theorem 1 is proved, following the steps summarized above.
The first part  (\S1--\S8) is ``purely combinatorial'', culminating in Propositions 2, 6 and 7 that are needed for the $SU(d)$ analysis.
The first part  is closely related to the proof of the discretized ring theorem in $\Bbb R$, presented in [B2];
the generalization to the higher dimensional setting necessitates reproducing  several technical  portions from that paper.

The counterpart in $\Bbb C$ of the main theorem from [B2] can be stated as follows.

\proclaim
{Proposition 2} Given $0<\sigma<2$ and $\kappa, \kappa'>0, \rho>0$, there are $\ve_0, \ve_0', \ve_1>0$ such that the following holds.

Let $A\subset\Bbb C\cap B(0, 1)$ satisfy

\roster
\item "{\rm (8.1)}" $N(A, \delta)= \delta^{-\sigma}$ \qquad ($\delta$ small enough).

\item "{\rm (8.2)}" $N(A\cap B(z, t), \delta)< t^\kappa N(A, \delta)$ if $\delta< t< \delta^{\ve_0} $ and $z\in\Bbb C$.

Let $\mu$ be a probability measure on $\Bbb C\cap B(0, 1)$ such that

\item "{\rm(8.3)}" $\mu \big(B(z, t)\big) < t^{\kappa'}$ if $ \delta < t<\delta^{\ve_0'}$ and $z\in\Bbb C$.

Let $z_1, z_2 \in\Bbb C$ satisfy

\item
"{(8.4)}" $\delta^{\ve_0'}< |z_1|\sim |z_2|< 1 $ and \Big|Im $\frac {z_1}{z_2}\Big| >\rho $

Then one of the following holds

\item "{\rm (8.5)}" $N(A+A, \delta)> \delta^{-\sigma-\ve_1}$.

\item "{\rm (8.6)}" $N(A+bA, \delta)> \delta^{-\sigma-\ve_1} $ for some $b\in \supp \mu$.

\item "{\rm (8.7)}" $N(A+z_1A, \delta)+N(A+z_2A, \delta)> \delta^{-\sigma-\ve_1}$.
\endroster
\endproclaim

\noindent
{\bf Acknowledgement.} The authors are grateful to Peter Sarnak for discussions related to \S10.
\bigskip

\centerline
{\bf Part 1: Generalizing the Discretized Ring Theorem}

The aim in what follows is to establish higher dimensional analogues of the discretized sum/product theory from [B2], in particular,  for subsets of $\Bbb C$ and $\Bbb
C^d$
\bigskip

\centerline
{\bf \S1. Basic Notation and Assumptions}

Let $\delta=2^{-m}$ ($m$ large) and $A\subset [0, 1]^d$ be a collection of $\delta$-separated points in $\Bbb R^d$ (alternatively we could take $A$ a union of $\delta$-intervals).
Assume
$$
|A|=\delta^{-\sigma}\tag {1.1}
$$
for some fixed $0<\sigma<d$ ($ | \ |$ denotes `cardinality' in (1.1) but may also be used for Lebesque measure if appropriate).

By size $\rho$ interval, we mean a $d$-dimensional box $\prod^d_{i=1} [a_i, a_i+\rho]$.

For $B\subset\Bbb R^d$ and $r>0$, we denote $N(B, r)$ the corresponding metrical entropy number.

We assume $A$ satisfies the following non-concentration property
$$
|A\cap I|< \rho^\kappa|A| \ \text { if } \ \delta<\rho< \delta^{\ve_0}  \ \text { and $I$ a size $\rho$ interval} \tag{1.2}
$$
for some $\kappa>0$ and $\ve_0 =\ve_0(\sigma)> 0$ small enough.

Let $\mu$ be a distribution on $\Cal L(\Bbb R^d, \Bbb R^d)$ (the space of linear maps on $\Bbb R^d$  satisfying certain assumptions to be specified (see Theorem 1 in \S7). 
Our aim is to show that for some $\ve_1>0$ (depending on the parameters), we have
$$
N(A+A, \delta)> \delta^{-\sigma-\ve_1} \text { or } \ N(A+bA, \delta)>\delta^{-\sigma -\ve_1} \text { for some $b\in \supp \mu$}.\tag {1.3}
$$
Let $T$ be a large constant (depending on the parameters and to be specified)
$$
m=Tm_1.\tag {1.4}
$$
We also consider a dyadic partition in intervals
$$
I_{n, k}= \prod^d_{i=1} [k_i 2^{-n}, (k_i+1) 2^{-n} [ \ \text { where } \ 0\leq k_i< 2^n \ \text { and } \ n\geq 0.
$$

We call $I_{n, k}$ a $2^{-n}$-interval.
\bigskip
\bigskip

\centerline
{\bf \S2. Initial Regularization of the Set}

We extract from $A$ a large subset with a `tree-structure'.
This construction is independent from assumptions (1.2), (1.3).

We introduce a subset $A_1\subset A$, a subset $\Cal S\subset\{1, \ldots, m_1\}$ and for $s\in\Cal S$ a dyadic integer 
$2\leq R_s< 2^{d(n_s-sT)}$, where $sT<n_s<(s+1) T-4$, with the following properties.

\roster
\item "{(2.1)}"  If $s\not\in \Cal S$ and $I$ is a $2^{-sT}$-interval, then there is at most
one $2^{-(s+1)T}$-interval $J\subset I$ such that $J\cap A_1\not=\phi$.
Take $R_s=1$ and $n_s =sT$ in this case.

\item "{(2.2)}"  If $s\in\Cal S$ and $I$ is a $2^{-sT}$-interval with $I\cap A_1\not=\phi$, then the number of $2^{-n_s}$-intervals
$J\subset I$ such that $J\cap A_1\not=\phi$ is $R_s$ and each such $J$-interval contains a single $2^{-(s+1)T}$-interval $J'$ intersecting
$A_1$.
Hence also
$$
R_s\sim N(A\cap I, 2^{-(s+1)T}).\tag 2.3
$$
Moreover, there is a pair of $2^{-n_s}$-intervals $J_1, J_2\subset I$ s.t. $J_1\cap A_1\not=\phi, J_2\cap A_1\not=\phi$ and
$$
2.2^{-n_s}< \text{dist\,} (J_1, J_2)< 10.2^{-n_s}.\tag 2.4
$$
\item"{(2.5)}" \ $ |A_1|=\prod_{s\in\Cal S} R_s> (cT^{-2})^{m_1} |A|>\delta^{-\sigma+
\frac{\log T}{T}}$.
\endroster

The construction is straightforward.

We start at the bottom of the tree, considering for each $2^{-(m_1 -1)T}$-interval $I$ s.t. $I\cap A\not= \phi$ the number of
$2^{-m_1 T}$-intervals $J\subset I$ with $J\cap A\not=\phi$.
If their number is less than $10^{3d}$,  fix one such $2^{-m_1 T}$-interval $J=J_I$.

If their number is larger than $10^{3d}$, introduce the largest integer (depending on $I$)
$$
(m_1-1) T<n< m_1T-4
$$ 
for which there is a pair of $2^{-n}$-intervals $J_1, J_2\subset I, J_1\cap A\not=\phi, J_2\cap A\not=\phi$ and
$$
2.2^{-n}< \text {dist\,} (J_1, J_2)< 10 \, {2^{-n}}.
$$
It is easily seen that from our definition of $n$
$$
N(A\cap I, 2^{-n})\sim N(A\cap I, 2^{-m_1 T}).
$$
Define $R$ the dyadic integer such that the number of $2^{-n}$-intervals $J\subset I, J\cap A\not= \phi$ is between $R$ and $2R$.
Thus
$$
N(A\cap I, 2^{-m_1 T})\sim R< 2^{d(n-(m_1-1)T)}.
$$
Obviously the integers $n$ and $R$ (depending on $I$) take at most $T$-values.
We may therefore clearly introduce a subset $A^{(m_1-1)}\subset A$,
$$
|A^{(m_1-1)}|> C T^{-2}|A|\tag 2.6
$$
satisfying one of the following alternatives:

\roster
\item "{(2.7)}"  If $I$ is an $2^{-(m_1-1)T}$-interval, there is at most one $2^{-m_1 T}$-interval
$J\subset I$ with $J\cap A_1\not=\phi$.   In this case, $m_1-1\not\in \Cal S$.

\item "{(2.8)}"  If $I$ is a $2^{-sT}$-interval with $I\cap A_1\not=\phi$, then (2.2) holds for some $n=n_{m_1-1}$ and $R=R_{m_1-1}$.  In this case, $m_1 -1\in\Cal S$.
\endroster

Next, repeat the construction for the set $A^{(m_1-1)}$ considering $2^{-(m_1-2)T}$-intervals $I$, $I\cap
A^{(m_1-1)}\not=\phi$ and $2^{-(m_1-1)T}$-subintervals $J\subset I$.
We obtain $A^{(m_1-2)}$ as the intersection of $A^{(m_1-1)}$ and a collection of $2^{-(m_1-1)T}$-intervals;
$$
|A^{(m_1-2)}|>cT^{-2}|A^{(m_1-1)}|\tag 2.9
$$
and $A^{(m_1-2)}$ satisfies either (2.1) or (2.2) with $s=m_1-2$, for some \hfill\break
$(m_1-2)T< n_{m_1-2} <(m_1-1)T-4$ and $R_{m_1-2}< 2^{d(n_{m_1-2}-(m_1-2)T)}$.

From construction
$$
A^{(m_1-2)}\cap I=A^{(m_1-1)}\cap I
$$
if $I$ is a $2^{-(m_1-1)T}$-interval intersecting $A^{(m_1-2)}$.
Hence properties (2.1), (2.2)  at level $m_1-1$ remain preserved.

The continuation of the process is clear.
We obtain
$$
A\supset A^{(m_1-1)}\supset A^{(m_1-2)}\supset\cdots\supset A^{(s)} \supset A^{(s-1)}\supset\cdots\supset
A^{(1)}
$$
where
$$
|A^{(s-1)}|> C T^{-2}|A^{(s)}|\tag 2.10
$$
and
$$
A^{(s-1)}\cap I=A^{(s)}\cap I
$$
if $I$ is a $2^{-sT}$-interval intersecting $A^{(s-1)}$.

Hence $A^{(s-1)}$ keeps the properties of $A^{(s)}$ at scales $2^{-n}$ for $n\geq sT$.

Let $A_1=A^{(1)}$.
Iteration of (2.10) gives (2.5).

Denote
$$
\Cal S =\{s_1< s_2< \cdots < s_t\}.
$$
For each $c\in \{0,1\}^t$, we will introduce an element $x_c\in A_1\subset A$.

Since $s_1\in \Cal S$, it follows from (2.2) that there is a pair of $2^{-n_{s_1}}$-intervals $I_0, I_1$
intersecting $A_1$ s.t.
$$
2.2^{-n_{s_1}}< \text {dist\,} (I_0, I_1)< 10.2^{-n_{s_1}}.\tag 2.11
$$
Denote $I_0'\subset I_0, I_1'\subset I_1$ the $2^{-s_2T}$-intervals intersecting $A_1$.
Again by (2.2), there are pairs $I_{0, 0}, I_{0, 1}\subset I_0'$ and $I_{1, 0}, I_{1,1 }\subset I_1'$ of
$2^{-n_{s_2}}$-intervals satisfying
$$
\align
2.2^{-n_{s_2}}&<\text {dist\,} (I_{0, 0}, I_{0, 1})< 10.2^{-n_{s_2}}\tag 2.12\\
2.2^{-n_{s_2}}&< \text {dist\,} (I_{1, 0}, I_{1, 1})<10.2^{-n_{s_2}}.\tag 2.12$'$
\endalign
$$
Continuing the construction, we obtain eventually $2^{-s_t T}$-intervals $I_c=I_{c_1, \ldots, c_t}$ for
$(c_1, \ldots, c_t)\in \{0, 1\}^t$ s.t. $I_c\cap A_1\not=\phi$ and with the property

If $c_1=c_1', \cdots, c_\tau =c_\tau', c_{\tau+1}\not= c_{\tau+1}'$, then
$$
2.2^{-n_{s_{\tau+1}}}< \text {dist\,}(I_c, I_{c'})< 10.2^{-n_{s_{\tau+1}}}.\tag 2.13
$$
Take
$$
x_c\in I_c\cap A_1 \text { for } c\in \{0, 1\}^t.\tag 2.14
$$
Fix some $\beta>0$ (to be specified) and define
$$
\Cal S_1 =\{1\leq s\leq m_1\big| (s+1)T-n_s >\beta T\}\tag 2.15
$$
(the porous levels). Obviously $\Cal S_1\supset \{1, \ldots, m_1\}\backslash \Cal S$.

Let $s_1< s_2< \cdots< s_{t_1}$ be an enumeration of the elements of $\Cal S_1$.

Let $\ve_2>0$ (to be specified) and assume first that
$$
t_1> \ve_2 m_1\tag 2.16
$$
which we refer to as the `porous case'.
The amplification in this situation can be performed exactly as in the $d=1$ case.
The argument is repeated in the next sections.

If $t_1\leq \ve_2m_1$, then $|\Cal S|>(1-\ve_2)m_1$ at most levels $s_\tau\in\Cal S$, (2.13) implies
$$
2.2^{-T} < 2^{s_\tau T} \dist (I_c, I_{c'}) < 10.2^{-(1-\beta)T}\tag 2.17
$$
if $c_1=c_1', \ldots, c_{\tau-1} = c_{\tau-1}'$ and $c_{\tau}\not= c_{\tau}'$.

This is the non-porous case and requires a different argument.
\bigskip

\bigskip
\centerline
{\bf \S3. The Porous Case (Bunching Together of Levels)}

Let $\Cal S=\{\Cal S_1<\cdots< \Cal S_{t_1}\}$ and assume (2.16).

We construct from $A_1$ a new system of sets $B_{k_1, \ldots, k_j}(j\leq t_1)$ with a `porosity property' at each level.

Denote $B_{k_1}, 1\leq k_1\leq K_1=\prod_{s\leq s_1} R_s$, the collection of non-empty intersections
$A_1\cap I$, where $I$ is a $2^{-n_{s_1}}$-interval.
By (2.1), (2.2), for each $B_{k_1}$ there is at most one $2^{-(s_1+1)T}$-interval $J$ such that
$\phi\not= A_1\cap J=B_{k_1}$.
Hence (reducing $|A_1|$ by a factor $c_d$) we may assume

\itemitem {(3.1)} \ diam$\,B_{k_1} \leq 2^{-(s_1+1)T}$

\noindent
and

\itemitem {(3.2)} \ dist\,$(B_{k_1}, B_{k_1'})> 2^{-(s_1+1-\beta)T}$ if $k_1\not= k_1'$.

Fixing $k_1$, let $B_{k_1k_2}, 1\leq k_2\leq \prod_{s_1<s\leq s_2} R_s= K_2$, be the collection of
non-empty intersections $B_{k_1}\cap I$ where $I$ is a $2^{-n_{s_2}}$-interval.
Again, each $B_{k_1, k_2}$ is contained in a single $2^{-(s_2+1)T}$-interval and

\itemitem{(3.3)} \ $B_{k_1}=\bigcup_{k_2\leq K_2} B_{k_1} 
$

\itemitem{(3.4)}  \ diam $B_{k_1k_2} < 2^{-(s_2+1)T}$

\noindent
and we may ensure

\itemitem {(3.5)} \ dist$(B_{k_1 k_2}, B_{k_1', k_2'})> 2^{-(s_2+1-\beta)T}$ for
$(k_1, k_2)\not= (k_1', k_2')$.

Continuing, we obtain a system $B_{k_1, k_2, \ldots , k_j} \ (1\leq k_i\leq K_i, j\leq
t_1)$ with the following properties

\itemitem{(3.6)} \ $\big|\bigcup\limits_{1\leq k_1\leq K_1} B_{k_1}\big|\geq C^{-t_1}|A_1| >C^{-m_1}|A_1|$

\itemitem{(3.7)} \ $ B_{k_1\ldots k_j}=\bigcup\limits_{1\leq k\leq K_{j+1}} B_{k_1\ldots k_jk}$

\itemitem{(3.8)} \  diam$\, B_{k_1\ldots k_j} \leq 2^{-T\beta}\lambda_j$

\itemitem{(3.9)} \ dist$ \,(B_{k_1, \ldots, k_j}, B_{k_1', \ldots, k_j'})>\lambda_j$ if $ \ \bark \ \not= \ \bark '$

\noindent
where we denoted
$$
\lambda_j=2^{-(s_j+ 1-\beta) T}.\tag 3.10
$$

Next we introduce a new system $C_{\ell_1\ldots\ell_s} \ (1\leq\ell_s\leq L_s$ and $ 1\leq s\leq t_2$).

Define
$$
M=[ 2^{T^{2/3}}]\tag 3.11
$$
and let
$$
10^3<D\in\Bbb Z_+\tag 3.12
$$
to be specified.

We consider the system $B_{k_1\ldots k_j}$ constructed above.

Starting from $t_1$, let $r_1 \in\Bbb Z_+$ be minimum s.t.
$$
K_{t_1} \ldots K_{t_1-r_1}> M^{r_1+1}.\tag 3.13
$$
We distinguish two cases.

If $r_1 \leq 10^3$, identify levels $t_1-r_1, \ldots, t_1$ to a single one, with branching
$$
L_{t_2} =K_{t_1}\ldots K_{t_1-r_1}> M^{r_1+1}.\tag 3.14
$$
The sets $C_{\ell_1\ldots\ell_{t_2}}$ are sets $B_{k_1\ldots k_{t_1}}$, hence satisfying

\itemitem {(3.15)} \ diam$\,C_{\ell_1\ldots \ell_{t_2}}< \sigma_{t_2}\mu_{t_2}$

\itemitem {(3.16)} \ dist$\,(C_{\ell_1\ldots, \ell_{t_2}}, C_{\ell_1', \ldots, \ell_{t_2}'})
>\mu_{t_2}$ if $ \overline\ell \not= \overline {\ell'}$

\noindent
where $\mu_{t_2} =\lambda_{t_1}$ and
$$
\sigma_{t_2}^{-1} = 2^{\beta T}> M^{D(r_1+1)}\tag 3.17
$$
assuming, if (3.11)
$$
T> 10^{10}\Big(\frac D\beta\Big)^3.\tag 3.18
$$
If $r_1> 10^3$, identify levels $t_1-r_1, \ldots, t_1- \frac{r_1}{100}+1$ to a single one, with branching
$$
L_{t_2} =K_{t_1-r_1}\cdots K_{t_1-\frac{r_1}{100}+ 1}>
\frac {M^{r_1+1}}{K_{r_1-\frac{r_1}{100}}\ldots K_{t_1}}> M^{r_1+1-(\frac{r_1}{100}+1)}> M^{\frac{99}{100}
r_1}\tag 3.19
$$
(from definition of $r_1$).
In order to ensure proper separation, reduce the sets
$B_{k_1, \ldots k_{t_1-\frac{r_1}{100}}}$ to their subset
$B_{k_1\ldots k_{t_1-\frac{r_1}{100}}, \underbrace{1\ldots, 1}_{\frac{r_1}{100}}} =C_{\ell_1\ldots\ell_{t_2}}$.

By (3.8), (3.9), the mutual distance between those sets is at least $\lambda_{t_1-\frac{r_1}{100}}\equiv
\mu_{t_2}$, while their diameter is at most $\lambda_{t_1}=\sigma_{t_2}\mu_{t_2}$ with
$$
\sigma_{t_2}^{-1} = 2^{\frac\beta{100}r_1 T}>M^{D(r_1+1)}\tag 3.20
$$
\big(by (3.18)\big).

We have reduced the size of $A_1$ by at most a factor $K_{t_1-\frac{r_1}{100}}\cdots K_{t_1}<
M^{\frac{r_1}{100}+1}< M^{r_1/90}$.

Next, repeat (if possible) the procedure, starting from level $t_1-r_1$ to obtain a level $t_1-r_1-r_2$, etc.

If at  some level $t'=t_1-r_1-r_2-\cdots$ we can not continue the process, it means that
$$
K_1\cdots K_{t'}\leq M^{t'+1}<
2^{T^{2/3} t_1}\leq 2^{m_1 T^{2/3}} <\Big(\frac 1\delta\Big)^{T^{-1/3}}.\tag 3.21
$$
Hence
$$
|A\cap I|\geq |B_{k_1\ldots k_{t'}}|\overset {(3.6)}\to> \delta^{T^{-1/3}}C^{-m_1}|A_1|\overset{(2.5)}\to>
\delta^{2T^{-1/3}}|A|\tag 3.22
$$
where $I$ is some $2^{-s_{t'}T}$-interval.

Recall the non-concentration assumption (1.2).
It follows from (3.22) that either $2^{-s_{t'}T}> \delta^{\ve_0}= 2^{-\ve_0m_1 T}$ or
$ \delta^{2T^{-1/3}}< 2^{-s_{t'}T\kappa}$.
Hence
$$
t'\leq s_{t'}< \max (\ve_0 m_1, \kappa^{-1} T^{-1/3}m_1)< \frac{\ve_2}{2} m_1\tag 3.23
$$
assuming
$$
\ve_0 <\frac{\ve_2}2 \text { and }  T^{-1/3}<\frac 12 \kappa\ve_2.\tag 3.24
$$
Since by (2.16), $t_1> \ve_2$, this will ensure that $t'<\frac 12 t_1$
and hence
$$
r_1+r_2+\ldots = t_1-t'> \frac 12 t_1>\frac 12\ve_2 m_1.\tag 3.25
$$
Consequently we replaced $B_{k_1\ldots k_j}(1\leq j\leq t_1)$ by a subtree
$C_{\ell_1\ldots\ell_s}(1\leq s\leq t_2)$ with the following properties

\itemitem {(3.26)} \ $C_{\ell_1\ldots\ell_{s-1}} =\bigcup\limits_{\ell\leq L_s} C_{\ell_1\ldots\ell_s\ell}\subset
A$

\itemitem {(3.27)} \ dist$\,(C_{\ell_1\ldots\ell_s}, C_{\ell_1',\ldots, \ell_s'})\geq \mu_s $ if $ (\ell_1,\ldots,
\ell_s)\not= (\ell_1', \ldots, \ell_s')$

\itemitem{(3.28)} \ diam$\, C_{\ell_1\ldots\ell_s}<\sigma_s\mu_s$, where $2^{-\beta T}>\sigma_s>\delta^{1/2}$

\itemitem {(3.29)} \ $\prod\limits_s(L_s \wedge \sigma_s^{-\frac 1D})> M^{\frac{99}{100}(r_1+r_2+\cdots)}
=M^{\gamma m_1}$

\noindent
where
$$
\gamma =\frac {99}{100 m_1}(r_1+r_2+\cdots)> \frac 13 \ve_2\tag 3.30
$$

\itemitem {(3.31)} \ $\big|C_\phi=\bigcup\limits_{\ell_1\leq L_1} C_{\ell_1}\big|> M^{-\frac
1{90}(r_1+r_2+\cdots)} c^{-m_1}|A_1|> M^{-\frac\gamma{80} m_1}|A|$.

The sets $C_{\overline\ell}$ are subsets of the discrete set $A$.

In what follows, it will be convenient to replace our discrete sets by unions of size-$\delta$ intervals,
defining
$$
A'=\Big\{x\in\Bbb R| \text { dist}(x, A)< \frac\delta 2\Big\}\tag 3.32
$$
and similarly, for $\overline\ell =(\ell_1, \ldots, \ell_s)$
$$
C_{\overline\ell}' =\Big\{x\in\Bbb R| \text {dist}(x, C_{\overline\ell})<\frac \delta 2\Big\}.\tag {3.33}
$$
Hence
$$
\align
|A'|&= \delta|A|=\delta^{1-\sigma}\tag 3.34\\
|C_{\overline\ell}'|&=\delta|C_{\overline\ell}|\tag 3.35
\endalign
$$
(using $| \ |$ to denote both measure and cardinality).

\bigskip

\centerline
{\bf \S4. Porous Case (Amplification)}

Assume $\mu $ a probability measure on $\Cal L(\Bbb R^d, \Bbb R^d)$ and $\kappa'>0$ satisfying the following conditions

\itemitem{(4.1)} $\Vert b\Vert \leq 1$ if $b\in \supp\mu$.

\itemitem {(4.2)} Given a unit vector $v\in\Bbb R^d$ and a vector $w\in\Bbb R^d$, we have
$$
\mu[b; |bv-w|<\rho]< c\rho^{\kappa'}
$$
for all $\delta<\rho<1$.

We denote by $\Bbb E=\Bbb E_b$ the $\mu$-expectation.
Our aim is to estimate from below
$$
|C_\phi'+ C_\phi'|+\Bbb E[|C_\phi'+ bC_\phi'|].
$$
Let us first show that we may replace the $b$-distribution as to ensure moreover the property

\itemitem{(4.3)} $b$ is invertible and $\Vert b^{-1}\Vert<3$ for $b\in\supp\mu$.

This is seen as follows.
Consider the map
$$
b\mapsto \frac 13 (2. 1\!\!1+b)=b'
$$
and the image distribution $\mu'$ of $\mu$.
Clearly (4.2) still holds and (4.3) is now satisfied.
Also
$$
|C_\phi'+b'C_\phi'|\leq |C_\phi'+C_\phi'+C_\phi'+C_\phi'+C_\phi'+ bC_\phi'|\lesssim \Big(\frac{|C_\phi+C_\phi|}{|C_\phi|}\Big)^6|C_\phi'+bC_\phi'|
$$
and hence
$$
\Bbb E'[|C_\phi'+b'C_\phi'|]\lesssim \Big(\frac{|C_\phi+C_\phi|}{|C_\phi|}\Big)^6 \Bbb E[|C_\phi'+bC_\phi'|].\tag {$*$}
$$
Thus it will suffice to establish a lower bound on $\Bbb E[|C_\phi'+bC_\phi'|]$ under assumptions (2.1)-(2.3).

Recalling (3.29), denote
$$
M_s =L_s\wedge \sigma_s^{-\frac 1D} \qquad (1\leq s\leq t_2)\tag 4.4
$$
satisfying
$$
\prod_{s=1}^{t_2} M_s> M^{\gamma m_1}.\tag 4.5
$$
We choose the parameter $D$ to satisfy
$$
D>\frac {10}{\kappa'}.\tag 4.6
$$

Starting from $s=1$, we have $C_\phi'=\bigcup_{\ell_1\leq L_1} C_{\ell_1}'$, where the $C_{\ell_1}'$ are contained in intervals of size $\sigma_1\mu_1$ and
separation $>\mu_1$.

Thus
$$
C_\phi'+bC_\phi'=\bigcup_{\Sb \ell_1\leq L_1\\ \ell_1'\leq L_1\endSb} (C_{\ell_1}'+bC_{\ell_1'}').
$$
Partition $[0, 1]^d$ in intervals $I_\alpha$ of size $\eta>\mu_1$ where $\eta$ is chosen such that
$$
\align
M_1&=\max_\alpha|\{ \ell_1\leq L_1; C_{\ell_1}' \cap I_\alpha \not=\phi \}|\\
&=|\{\ell_1 \leq L_1; C_{\ell_1}'\cap I_0 \not=\phi\}|.\tag 4.7
\endalign
$$
Denote for each $\alpha$
$$
\Cal E_\alpha =\{\ell_1 \leq L_1; C_{\ell_1}'\cap I_\alpha\not= \phi\}
$$
and
$$
D_\alpha =\bigcup_{\ell_1\in \Cal E_\alpha} C_{\ell_1}' \subset I_\alpha +B(0, \sigma_1\mu_1).
$$
Then
$$
C_\phi' +bC_\phi'\subset \bigcup_\alpha (D_\alpha+bD_0)
$$
and from the preceding
$$
|C_\phi'+bC_\phi'|> c
\sum_\alpha |D_\alpha+bD_0|.\tag 4.8
$$
Fixing $\alpha$, we have
$$
D_\alpha+bD_0=\bigcup_{\Sb \ell_1\in \Cal E_\alpha\\ \ell_1' \in\Cal E_0\endSb} (C_{\ell_1}' +bC_{\ell_1'}')\tag 4.9
$$
and certainly
$$
|D_\alpha+bD_0|> c\max_{\ell_1' \in\Cal E_0}\sum_{\ell_1\in\Cal E_\alpha} |C_{\ell_1}'+bC_{\ell_1'}'|>\frac c{M_1} \sum_{\ell_1\in \Cal E_\alpha, \ell_1'\in \Cal
E_0} |C_{\ell_1}'+b C_{\ell_1'}'|.\tag 4.10
$$
On the other hand, fixing points $\xi_\ell \in C_\ell$, we have
$$
|\xi_\ell -\xi_{\ell'}|> \mu_1 \ \text { for } \  \ell\not=\ell'
$$
and
$$
C_\ell'+bC_{\ell'}' \subset \xi_\ell +b\xi_{\ell'}+B(0, 2\sigma_1\mu_1).\tag 4.11
$$
Therefore the sets in (4.9) will be mutually disjoint if $$
\min_{\Sb k, \ell\in\Cal E_r\\ k', \ell'\in\Cal E_0\endSb} |(\xi_k+b\xi_{k'})-(\xi_\ell+b\xi_{\ell'})|> 4\sigma_1\mu_1.\tag 4.12
$$

Here condition (4.2) comes into play.

Note that if $k'=\ell'$, $|\xi_k-\xi_\ell|>\mu_1\gg \sigma _1\mu_1$, so that we may assume $k'\not= \ell'$.

Fix $k, \ell \in\Cal E_r, k'\not=\ell '\in \Cal E_0$ and denote $v=\frac{\xi_{k'}-\xi_{\ell'}}{|\xi_{k'}-\xi_{\ell'}|}$.
It follows from (4.3) that the $\mu$-measure of the $b$'s for which
$$
\Big|bv+\frac{\xi_k-\xi_\ell}{|\xi_{k'}-\xi_{\ell'}|}\Big| < 4\sigma_1
$$
is at most $C\sigma_1^\kappa$.
Hence (4.12) holds for $b$  outside a set of $\mu$-measure at most
$$
C.|\Cal E_0|^2. |\Cal E_r|^2 \sigma_1^\kappa\leq CM_1^4 \sigma_1^\kappa <\sigma_1^{\kappa/2}\tag 4.13
$$
by (4.4), (4.5). For such $b$, we get
$$
|D_\alpha+bD_0|=\sum_{\Sb \ell_1\in \Cal E_\alpha\\ \ell_1'\in\Cal E_0\endSb} |C_{\ell_1}'+bC'_{\ell_1'}|.\tag 4.14
$$
From (4.10), (4.14), it follows that
$$
|D_\alpha+bD_0|>\sum_{\Sb \ell_1\in \Cal E_\alpha\\ \ell_1'\in\Cal E_0\endSb} \vp_{\ell_1, \ell_1'}(b)|C_{\ell_1}'+ bC_{\ell_1'}'|\tag 4.15
$$
where $\vp_{\ell_1, \ell_1'}$ takes values in $\{1, \frac c{M_1}\}$ and by (4.13)
$$
\mu[\vp_{\ell_1, \ell_1'}\not= 1]< \sigma_1^{\kappa/2}.\tag 4.16
$$
Summing over $\alpha$. (4.8), (4.15) imply
$$
|C_\phi'+bC_\phi'|> c\sum_{\ell_1\leq L_1, \ell_1'\in \Cal E_0}\vp_{\ell_1, \ell_1'} (b)|C_{\ell_1}'+bC_{\ell_1'}'|\tag 4.17
$$
with $\vp_{\ell_1, \ell_1'}$ as above.

Next, restrict $(\ell_1, \ell_1')$ to the set $\{1\leq \ell_1\leq L_1\}\times (\{1\leq \ell_1'\leq L_1\}\backslash \Cal E_0)$
and repeat the construction.
Note that if $(\ell_1, \ell_1')$ is restricted to a product set $\Cal F\times\Cal F'$, we partition in interval $I_\alpha$ of size $\eta>\mu_1$ chosen such that
$$
M_1= \max_\alpha|\{\ell_1\in\Cal F; C_{\ell_1}' \cap I_\alpha\not=\phi\}| \vee \max_\alpha |\{\ell_1'\in\Cal F'; C_{\ell_1'}'\cap I_\alpha\not= \phi\}|
$$
and obtained as either a set $\Cal E_0 =\{\ell_1\in\Cal F; C_{\ell_1}' \cap I_0\not= \phi\}$ or $\Cal E_0=\{\ell_1' \in\Cal F'; C_{\ell_1'}'\cap I_0\not=\phi\}$.
Because of assumption (4.3), both cases may be treated similarly.

Exhausting the set $\{1\leq \ell_1\leq L_1\}\times\{1\leq \ell_1'\leq L_1\}$ in $\sim\frac{L_1}{M_1}$ steps, we get
$$
|C_\phi'+bC_\phi'|> c\frac{M_1}{L_1} \sum_{\ell_1\leq L_1, \ell_1'\leq L_1} \vp_{\ell_1, \ell_1'}(b)|C_{\ell_1}'+b C'_{ \ell_1'}|\tag 4.18
$$
with $\vp_{\ell_1, \ell_1'}$ satisfying (4.13).  Equivalently
$$
|C_\phi'+ bC_\phi'|>\sum_{\ell_1\leq L_1, \ell_1'\leq L_1} \psi_{\ell_1, \ell_1'} (b)|C_{\ell_1}'+bC'_{\ell_1'}|\tag 4.19
$$
with $\psi_{\ell_1, \ell_1'}$ taking values in $\{c\frac{M_1}{L_1}, c\frac 1{L_1}\}$ and
$$
\mu\Big[\psi_{\ell_1, \ell_1'}\not= c\frac{M_1}{L_1}\Big]< \sigma_1 ^{\kappa'/2}.\tag 4.20
$$
Therefore
$$
\Bbb E[\log\psi_{\ell_1, \ell_1'}]> (1-\sigma_1^{\kappa'/2})\Big(\log c\frac{M_1}{L_1}\Big)
+\sigma_1^{\kappa'/2}\Big(\log \frac c{L_1}\Big) > \log \frac{\sqrt{M_1}}{L_1}.\tag 4.21
$$
We assume here $\sigma_1^{\kappa'/2}< \frac 12$, which will be fulfilled if
$$
T>\frac{10}{\beta\kappa'}\tag 4.22
$$
(since $\sigma_s< 2^{-\beta T}$).

Repeat with the sets
$$
 C_{\ell_1}'+bC'_{\ell_1'} =\bigcup_{\ell_2,\ell_2'\lesssim L_2} (C_{\ell_1, \ell_2}'+bC'_{\ell_1', \ell_2'})
$$
to obtain
$$
|C'_{\ell_1} +b C'_{\ell_1'}|> \sum_{\ell_2, \ell_2'}\psi_{\ell_{1}, \ell_2; \ell_1',
\ell_2'}(b)|C'_{\ell_1, \ell_2}+ bC'_{\ell_1', \ell_2'}|\tag 4.23
$$
where again
$$
\Bbb E_x[\log \psi_{\ell_1, \ell_2; \ell_1', \ell_2'}]>\log\Big(\frac 1{L_2}M_2^{1/2}\Big).
\tag 4.24
$$
Iteration clearly provides the following minoration
$$
\align
|C_\phi'+ bC_\phi'|&>\sum_{\Sb \ell_1, \ldots, \ell_{t_2}\\ \ell_1', \ldots, \ell_{t_2}'\endSb}
\psi_{\ell_1, \ell_1'}(b) \cdots\psi_{\ell_1, \ldots, \ell_{t_2}; \ell_1', \ldots, \ell_{t_2}' }
(b)| C'_{\ell_1, \ldots, \ell_{t_2}}+bC'_{\ell_1',\ldots, \ell_{t_2}'}|\\
& \geq \sum_{\overline\ell,\overline \ell'}\psi_{\ell_1, \ell_1'}\cdots \psi_{\overline \ell; \overline
\ell}|C'_{\ell_1\cdots\ell_{t_2}|}\tag 4.25
\endalign
$$
where by (4.21), (4.24),  etc.
$$
\align
\Bbb E_x[\psi_{\ell_1, \ell_1'}\psi_{\ell_1, \ell_2; \ell_1', \ell_2'\ldots}]
&\geq e^{\Bbb E[\log\psi_{\ell_1, \ell_1'}]+\Bbb E[\log\psi_{\ell_1, \ell_2; \ell_1',
\ell_2'}]+\cdot}\\
&> \prod^{t_2}_{s=1} \frac{M_s^{1/2}}{L_s}.\tag 4.26
\endalign
$$
Therefore
$$
\align
\Bbb E[|C_\phi'+bC_\phi'|]& >\Big(\prod^{t_2}_{s=1} M_s^{1/2}\Big) \sum_{\ell_1, \ldots, \ell_{t_2}}|C_{\ell_1, \ldots, \ell_{t_2}}'|\\
&>M^{\frac\gamma 2 m_1} |C_\phi'|.\tag 4.27
\endalign
$$
Recall $(*)$, it follows that if $\mu$ satisfies (4.1) and (4.2)
$$
|C_\phi'+ C_\phi'|+\Bbb E[|C_\phi'+bC_\phi'|]> M^{\frac\gamma{14} m_1}|C_\phi'|.\tag 4.28
$$
Therefore, by (3.31), (3.30)
$$
\align
N(A+A, \delta)+\Bbb E[N(A+bA, \delta)]&> M^{\frac\gamma{14} m_1}|C_\phi|\\
&>M^{\frac\gamma{20} m_1} |A|\\
&> \delta^{-\frac 1{60} \ve_2 T^{-1/3}}|A|.\tag 4.29
\endalign
$$
\bigskip

\centerline
{\bf \S5. The Non-Porous Case (1)}

Assume $|\Cal S_1|\leq \ve_2 m_1$ and denote
$$
\Cal S'= \{1, \ldots, m\}\backslash\Cal S_1\subset\Cal S\tag 5.1
$$
satisfying
$$
t'= |\Cal S'|> (1-\ve_2) m_1\tag 5.2
$$
and (2.17) if $\tau\in\frak Z'= \{1\leq \tau\leq t; s_\tau\in\Cal S'\}\quad (t=|\Cal S|)$.

For notational simplicity, we identify $\Cal S'$ with $\{ 1, \ldots, m_1\}$.
Recalling the system $\{x_c\}_{c\in\{ 0, 1\}^t}$ of points from $A$ introduced in (2.14)
and (2.17), our starting point is a system $\{x_c\}_{c\in\{0, 1\}^{m_1}}\subset A$ such that if $c_1=c_1', 
\ldots, c_s = c_s', c_{s+1}'\not= c_{s+1}'$, 
then 
$$
2.2^{-(s+1)T} < |x_c-x_{c'}|< 10.2 ^{-(s+1-\beta)T}.\tag 5.3
$$
Denote $J=[\frac 1{10d}2^{(1-\beta)T}]$ and $k=k(T)\in\Bbb Z_+$.

We construct in the $k$-fold sumset $kA$ of $A$, for $s=1, \ldots, m$, points
$$
(y_{j_1, \ldots, j_s})_{\Sb 1\leq j_1\leq J\\ \vdots\\ 1\leq j_s\leq J\endSb}\tag 5.4
$$
with the following properties

\roster
\item "{(5.5)}" $|y_{j_1\ldots, j_s j_{s+1}} -y_{j_1\ldots j_s}|<\frac 1{10d} 2^{-sT}$

\item "{(5.6)}" $2^{-(s+1)T} < |y_{j_1\ldots j_s 1}- y_{j_1, \ldots j_s}|<\frac 1{10dJ} 2^{-sT}$

\item "{(5.7)}" For fixed $j_1, \ldots, j_s$, we have (setting $y_\phi=0$ for $s=0$)
$$
|(y_{j_1\ldots j_{s+1}} - y_{j_1}\ldots j_s)-j_{s+1}(y_{j_1\ldots j_s 1} -y_{j_1\ldots j_s})|< 4^{-T} 2^{-sT}.
$$
\endroster

Thus (5.7) means that the points $(y_{j_1}\ldots _{j_{s} j_{s+1}})_{1\leq j_{s+1}\leq J}$ lie approximatively on some 1-dim line segment emanating from
$y_{j_1}\ldots _{j_s}$ and length between $2^{-sT-\beta T}$ and $2^{-sT}$

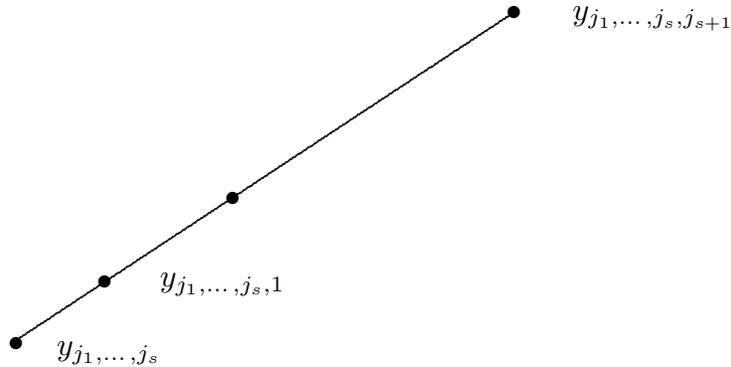
\input fig1.tex

\bigskip

\noindent
The construction is done as follows

The points $y_{j_1\ldots j_s}$ will be of the form
$$
y_{j_1\ldots j_s } = x_{c^{(1)} }+ \cdots+ x_{c^{(k)}}\tag 5.8
$$
for certain $c^{(1)}, \ldots, c^{(k)} \in  \{0, 1\}^{m_1}$ such that $c_{s'}^{(1)} =\cdots= c_{s'}^{(k)} =0$ for $s' >s$.

Denote $d^{(1)}, \ldots, d^{(k)} \in \{0, 1\}^{m_1}$ by
$$
\left\{
\aligned
&d_{s'}^{(\ell)} = \,  c_{s'}^{(\ell)} \  \text { if } \ s'\leq s\cr
&d_{s+1}^{(\ell)} = \, 1\cr
&d_{s'}^{(\ell)} =  \, 0  \ \text { if } \ s'>s+1.
\endaligned\right.
$$

From (5.3), for each $\ell =1, \ldots, k$
$$
2^{-(s+1)T}< |x_{c^{(\ell)}} - x_{d^{(\ell)}}|<\frac 1{10dJ} 2^{-sT}.\tag 5.9
$$
Taking $k=k(T)$ sufficiently large $(\log k(T)\sim T)$, we may clearly specify a subset $\Cal L=\{\ell_1< \cdots< \ell_J\}\subset \{ 1, \ldots,  k\}$ 
such that for all $\ell\in\Cal L$
$$
|x_{c^{(\ell)}} - x_{d^{(\ell)}} -\eta 2^{-(s+1)T} v|< 2^{-sT-3T}\tag 5.10
$$
for some unit vector $v\in\Bbb R^d$ and some $1<\eta< 2^{\beta T}$.

Define
$$
y_{j_1\ldots j_sj}=\sum_{\ell\not\in \{\ell_1, \ldots, \ell_j\}} x_{c^{(\ell)}} 
+\sum_{\ell\in\{\ell_1, \ldots, \ell_j\}} x_{d^{(\ell)}}.\tag 5.11
$$
Hence (5.5), (5.6) follow from (5.8), (5.11) and (5.9).

By (5.10)
$$
\align
y_{j_1\ldots j_{s}j}-y_{j_1\ldots j_s}&= (x_{d^{(\ell_1)}}-x_{c^{(\ell_1)}})+\cdots+ (x_{d^{(\ell_j)}}-x_{c^{(\ell_j)}})\\
&= j\eta 2^{-(s+1)T}v+ 0(j 2^{-sT-3T})
\endalign
$$
implying (5.7).

\bigskip
\centerline
{\bf \S6. Non-Porous case (2)}

We make the following further assumption on the distribution $\mu$ of $b\in \Cal L(\Bbb R^d, \Bbb R^d)$.

\itemitem{(6.1)} {\sl There is a function $\theta (\rho)\to 0$ for $\rho\to 0$ such that if $v, w$ are unit vectors in $\Bbb R^d$ then
$$
\mu[b|
\langle bv,w\rangle|<\rho]< \theta(\rho).
$$}

Strictly speaking, all we require is that for some $\rho>0$
$$
\max_{|v|=1=|w|}\mu[b| \ |\langle bv, w\rangle| < \rho]<\theta\tag 6.2
$$
with $\theta$ sufficiently small.

Fix a constant
$$
\ve_3>0\tag 6.3
$$
(to be specified).

Property (6.2) will be used in the following

\proclaim
{Lemma 6.4}
Let $v_1, \ldots, v_d\in \Bbb R^d$ be unit vectors and $\rho^{-1} 2^{-T}< \eta_1, \ldots, \eta_d <\frac 1J$.
Define
$$
P_i=\{j\eta_i v_i; 1\leq j\leq J\} \ \text { for } \ i=1, \ldots, d.
$$
Then
$$
\Bbb E[\min_{x\not= y\in b_1 P_1+\cdots+ b_d P_d} |x-y|< 2^{-T}]<\ve_3\tag 6.5
$$
where $\Bbb E$ refers to the $d$-fold product measure $\mu\otimes \cdots\otimes \mu =\mu^{(d)}$.
\endproclaim

\noindent
{\bf Proof.}

Write
$$
b_1P+\cdots+b_dP =\Big\{\sum^d_{i=1} j_i\eta_i b_iv_i; 1\leq j_i\leq J\Big\}
$$
We need to ensure that for $|j_1|+\cdots+ |j_d|= 0$
$$
|j_1\eta_1 b_1 v_1+\cdots+ j_d\eta_db_d v_d|> 2^{-T}.\tag 6.6
$$
It suffices to impose the conditions
$$
\left\{
\aligned
&|b_1v_1|>\rho\\
&\dist (b_2v_2, [b_1v_1])>\rho\\
& \qquad \vdots\\
&\dist (b_dv_d, [b_1v_1, \ldots, b_{d-1} v_{d-1}])>\rho.
\endaligned
\right.
$$

By (6.2), clearly
$$
\mu^{(d)}[(b_1, \ldots, b_d)| (6.7) \text { fails}]< d\theta <\ve_3
$$
for appropriate $\theta$.
This proves the Lemma.\qed

Let $\{y_{j_1\ldots j_s}\}$ be the systems obtained in \S5 satisfying (5.5)-(5.7).
We denote $Y=\{y_{j_1}\ldots _{j_{m_1}}\}$ and $Y_{j_1\ldots j_s} =\{y_{j_1\ldots j_s j_{s+1}\ldots j_{m_1}}\}$.
By (5.5)
$$
Y_{j_1\ldots j_s} \subset B(y_{j_1\ldots j_s}, \frac 1{5d}, 2^{-sT}).\tag 6.8
$$
Denoting $\Bbb E$ the expectation wrt $\mu^{(d)}$, we establish a lower bound on
$$
\Bbb E [N(b_1Y+\cdots+ b_dY, \delta)]
$$
following an argument similar to that used in \S4.

From (5.7), (6.5), the expectation for the points 
$b_1 y_{j_1^{(1)}}+\cdots+ b_d y_{j^{(d)}_1}$ to be at least
$2^{-T}$ apart is $>1-\ve_3$ (we use here that $|y_j-jy_1|< 4^{-T})$.
Hence, by (6.8)
$$
N(b_1Y+\cdots +b_dY, \delta)< \vp(b_1, \ldots, b_d) \sum_{1\leq j_1^{(i)}\leq J} N(b_1 Y_{j_1^{(1)}}+\ldots + b_d Y_{j_1^{(d)}},  \delta)\tag 6.9
$$
where $\vp(b)\in \{1, J^{-d}\}, \Bbb E[\vp=1]> 1-\ve_3$.
Thus
$$
\Bbb E[\log \vp]> \ve_3 \log J^{-d} =-\log J^{d{\ve_3}}.\tag 6.10
$$
Similarly
$$
\align
&N(b_1Y_{j_1^{(1)}}+\cdots+ b_d Y_{j_1^{(d)}}, \delta)>\\
&\vp_{j_1^{(1)} \ldots j_1^{(d)}} (b) \sum_{1\leq j_2^{(i)}\leq J} N(b_1 Y_{j_1^{(1)}j_2^{(1)}} +\cdots+ b_dY_{j_1^{(d)}j_2^{(d)}}, \delta)\tag 6.11
\endalign
$$
with
$$
\Bbb E[\log \vp_{j_1^{(1)}\ldots j_1^{(d)}}]> -\log J^{d\ve_3}.\tag 6.12
$$
Iterating, we obtain
$$
N(b_1Y_1+\cdots+ b_dY_d, \delta)>\sum_{J^{(1)}, \ldots, J^{(d)}} \vp(b) \vp_{j_1^{(1)}\ldots j_1^{(d)}}(b) \quad \vp_{j_1^{(1)}j_2^{(1)} \ldots
j_1^{(d)}j_2^{(d)}}(b)\ldots\tag 6.13
$$
where $J^{(i)}=(j_1^{(i)}, \ldots, j_{m_1}^{(i)})$, and the expectation of the summands in (6.13) is at least
$$
e^{\Bbb E[\log \vp]+\Bbb E[\log \vp_{j_1^{(1)}\ldots j_1^{(d)}}]+\ldots} > J^{-\ve_3 dm_1}.\tag 6.14
$$
Hence
$$
\Bbb E[N(b_1Y+\cdots+ b_dY, \delta)]> J^{dm_1-\ve_3 dm_1}>[c2^{(1-\beta)T}]^{(1-\ve_3)dm_1}.\tag 6.15
$$
Recalling that in \S5 we identified $\Cal S'$ with $\{1, \ldots, m_1\}$ and (5.2), (6.5) gives
$$
\align
\Bbb E[N(b_1 \tilde A+\cdots+b_d\tilde A, \delta)]&> [c2^{(1-\beta)T}]^{(1-\ve_2)(1-\ve_3) dm_1}\\
&> \delta^{-(1-\beta)(1-\ve_2)(1-\ve_3)(1-\frac c T)d}\tag 6.16
\endalign
$$
where $\tilde A$ is the sum set $kA\supset Y$.

Since again from the sumset inequalities
$$
N(b_1\tilde A+\cdots+ b_d\tilde A, \delta)\geq \Big[ \prod^d_{i=1} N(A+b_iA, \delta)\Big]^k . \frac 1{|A|^{kd-1}}\tag 6.17
$$
it follows from (6.16) that for some $b\in \supp\mu$
$$
\align
N(A+bA, \delta)&> |A|^{1-\frac 1{kd}}\delta^{-\frac 1k (1-\beta)(1-\ve_2) (1-\ve_3)(1-\frac cT)}\\
&>\delta^{\frac\sigma{kd}-\frac 1k (1-\beta)(1-\ve_2)(1-\ve_3)(1-\frac cT)}|A|.\tag 6.18
\endalign
$$
Recalling the hypothesis in Lemma 6.4, we assume
$$
2^{-\beta T} <\rho.\tag 6.19
$$

\bigskip
\centerline
{\bf \S7. Summary}

Recall the parameters

{$\sigma$} \ \qquad (1.1)

$\kappa, \ve_0$ \quad (1.2)

$T$ \qquad  (1.4)

$\ve_1$ \qquad  (1.3)

$\ve_2$ \qquad  (2.16)

$\kappa'$ \qquad (4.2)

$k$ \ \qquad \S5

$\theta, \rho$ \quad (6.2)

$\ve_3$ \quad \ (6.3)

\noindent
and the conditions (3.18), (3.24), (4.6), (4.22), (6.19).

In view of (6.18), take
$$
\beta=\ve_2=\ve_3 =10^{-2} \Big( 1 - \frac \sigma d\Big)\tag 7.1
$$
(recall that $0<\sigma <d$) and assume in (1.2)
$$
\ve_0< \frac 1{200} \Big(1-\frac\sigma d\Big).\tag 7.2
$$
In (6.2), $\rho$ depends on $\ve_3$, hence on $1-\frac \sigma d$.
From (3.18), (3.24), (4.6), (6.19), we impose on $T$ the condition
$$
T> 10^{13} (\kappa'\beta)^{-3} +10 (\kappa\ve_2)^{-3} +\frac 1\beta \log \frac 1\rho\tag 7.3
$$
taking
$$
T\sim 10^{20} \Big( 1-\frac \sigma d\Big)^{-3} \Big(\kappa^{-3}+(\kappa')^{-3}+\log \frac 1\rho\Big).\tag 7.4
$$
Recall also that $k=k(T)$, $\log k\sim T$.

From (4.29), (6.18), it follows
$$
\align
N(A+A, \delta)+ N(A+bA, \delta)&
>\min (\delta^{-\frac 1{60} \ve_2T^{-1/3}}, \delta^{\frac \sigma{kd} -\frac 1k(1-\beta)(1-\ve_2)(1-\ve_3) (1-\frac cT)}).|A|\\
&> \delta^{- \tau(\sigma, \kappa, \kappa', \rho)}|A|\tag 7.5
\endalign
$$
for some $b\in \supp\mu$.

We proved the following

\proclaim
{Proposition 1}

Let $A\subset [0, 1]^d$ and $N(A, \delta)=\delta^{-\sigma}$ for some fixed $0<\sigma< d$.

Assume for some $\kappa>0$

\roster
\item
"{\rm (7.6)}"  $N(A\cap I, \delta)< \delta_1^\kappa N(A, \delta) \text { if } \delta<\delta_1<\delta^{\ve_0}$ and $I\subset\Bbb R^d$ a size $\delta_1$-interval
\big(where $\ve_0=\frac 1{200} (1-\frac \sigma d)\big)$. 

Let further $\mu$ be a probability measure on $\Cal L(\Bbb R^d, \Bbb R^d)$ satisfying

\item "{\rm (7.7)}" \ $\Vert b\Vert\leq 1$ for $b\in \supp\mu$
\item "{\rm (7.8)}" \ There is $\kappa'>0$ such that
$$
\max_{\Sb v, w\in\Bbb R^d\\ |v|=1\endSb} \mu[|bv-w| < \delta_1]> c\delta_1^{\kappa'}\  \text { for } \ \delta <\delta_1<1
$$

\item "{\rm (7.9)}" \ There is $\rho>0$ such that
$$
\max _{|v|=1=|w|} \mu[|\langle bv, w\rangle| <\rho]<\frac 1{100d} \Big(1-\frac\sigma d\Big).
$$
\endroster

Then there is some $b\in \supp\mu$ such that
$$
N(A+A, \delta)+ N(A+bA, \delta)> \delta^{-\sigma-\tau}\tag 7.10
$$
where $\tau=\tau (\sigma,  \kappa, \kappa', \rho)>0$.
\endproclaim

Assume $u_1, \ldots, u_r\in\Cal L (\Bbb R^d, \Bbb R^d)$ such that  $\Vert u_i\Vert\leq 1$ and for any unit vectors $v, w\in\Bbb R^d$
$$
\max_{1\leq s\leq r}|\langle u_s v, w\rangle|>\rho\tag 7.11
$$
(note that $r$ may be restricted to some $r(d, \rho)$).

Let $J=J(d, \sigma)\in\Bbb Z_+$ be sufficiently large and let $\mu_1$ be the normalized image measure on $\Cal L(\Bbb R^d, \Bbb R^d)$ under the map
$$
(j_1, \ldots, j_r) \mapsto b=j_1 u_1+\cdots+ j_ru_r\qquad (1\leq j_s\leq J).
$$
If $v, w\in\Bbb R^d, |v|=1 = |w|$ and $t$ a scalar, we have, assuming $|\langle u_1 v, w\rangle|>\rho$,
$$
\align
\mu_1\Big[|t+bv.w|<\frac \rho2\Big] & = J^{-r} \sum_{j_2, \ldots, j_r} \Big|\Big\{ j_1\leq J; |t+j_1(u_1 v.w)+j_2(u_2v.w)+\cdots+ j_r(u_r v.w)|<\frac \rho
2\Big\}\Big|\\
&< \frac 1J.
\endalign
$$
Let $\mu_0$ on $\Cal L(\Bbb R^d, \Bbb R^d)$ satisfy (7.7), (7.8).
Then the image measure $\mu$ of $\mu_0\otimes\mu_1$ under the map $(b, b')\mapsto b+b'$ will clearly satisfy both (7.8), (7.9).

Thus one has

\proclaim
{Proposition 1$'$}

Let $A\subset [0, 1]^d$ and $N(A, \delta)=\delta^{-\sigma}$ for some $0<\sigma< d$.

Assume $A$ satisfies the non-concentration property (7.6).

Let $\mu$ be a probability measure on $\Cal L(\Bbb R^d, \Bbb R^d)$ satisfying (7.7), (7.8).

Let further $u_1, \ldots, u_r \in\Cal L(\Bbb R^d, \Bbb R^d), \Vert u_s\Vert\leq 1$ and $\rho>0$ such that
$$
\min_{|u|=1=|v|} \max_s |\langle u_s v, w\rangle|>\rho.\tag 7.12
$$
Then either
$$
N(A+A, \delta)> \delta^{-\sigma-\tau}
$$
or
$$
N(A+bA, \delta)>\delta^{-\sigma-\tau} \text { for some $b\in \supp\mu$}
$$
or
$$
N(A+u_sA, \delta)>\delta^{-\sigma-\tau} \text  { for some $s=1, \ldots, r$}
$$
where $\tau =\tau (\sigma, \kappa, \kappa', \rho)>0$.
\endproclaim
\bigskip\bigskip

\centerline
{\bf \S8. Discretized  Ring Theorem in $\Bbb C$}

Using Proposition 1$'$, we prove the following

\proclaim
{Proposition 2} Given $0<\sigma<2$ and $\kappa, \kappa'>0, \rho>0$, there are $\ve_0, \ve_0', \ve_1>0$ such that the following holds.

Let $A\subset\Bbb C\cap B(0, 1)$ satisfy

\roster
\item "{\rm (8.1)}" $N(A, \delta)= \delta^{-\sigma}$ \qquad ($\delta$ small enough).

\item "{\rm (8.2)}" $N(A\cap B(z, t), \delta)< t^\kappa N(A, \delta)$ if $\delta< t< \delta^{\ve_0} $ and $z\in\Bbb C$.

Let $\mu$ be a probability measure on $\Bbb C\cap B(0, 1)$ such that

\item "{\rm(8.3)}" $\mu \big(B(z, t)\big) < t^{\kappa'}$ if $ \delta < t<\delta^{\ve_0'}$ and $z\in\Bbb C$.

Let $z_1, z_2 \in\Bbb C$ satisfy

\item
"{(8.4)}" $\delta^{\ve_0'}< |z_1|\sim |z_2|< 1 $ and \Big|Im $\frac {z_1}{z_2}\Big| >\rho $

Then one of the following holds

\item "{\rm (8.5)}" $N(A+A, \delta)> \delta^{-\sigma-\ve_1}$.

\item "{\rm (8.6)}" $N(A+bA, \delta)> \delta^{-\sigma-\ve_1} $ for some $b\in \supp \mu$.

\item "{\rm (8.7)}" $N(A+z_1A, \delta)+N(A+z_2A, \delta)> \delta^{-\sigma-\ve_1}$.
\endroster
\endproclaim

This may be seen as the extension to $\Bbb C$ of the main result from [B2].

\noindent
{\bf Proof.}

We identify $\Bbb C$ with $\Bbb R^2$ viewing complex multiplication by $z=x+iy$ as
$$
\pmatrix x&-y\\y&x\endpmatrix \in\Cal L(\Bbb R^2, \Bbb R^2).
$$
Condition (8.3) has to be upgraded to (7.8) (i.e. removing the restriction $t<\delta^{\ve_0'}$).
We proceed as follows.

Define
$$
t_0=\text { inf } \{ t>\delta;\max_z \mu\big(B(z, t)\big)\geq t^{\frac 12\kappa'}\}
$$
obtained for $z=b$ say.
It follows from (8.3) that $t_0\geq \delta^{\ve_0'}$.
Denote
$$
\mu_1 =\frac {^{\mu|}B(b, t_0)}{\mu(B( b, t_0))}.
$$
From the definition of $t_0$, it follows that if $\frac \delta{t_0}< t\leq 1$ and $z\in\Bbb C$
$$
\mu_1 \big(B(z, tt_0)\big)<\frac {(t t_0)^{\frac 12 \kappa'}}{\mu \big(B(b, t_0)\big)} = t^{\frac 12\kappa'}.\tag 8.8
$$
In particular, there are elements $b', b''\in B(b, t_0)\cap \supp\mu$ such that $|b'-b''|> ct_0$.

Let $\mu_2$ be the image measure of $\mu_1$ under the map $z\mapsto \frac{z-b'}{b''-b'}$.
Clearly $\supp\mu_2 \subset B(0, C)$ and from (8.8)
$$
\sup_z \mu_2 \big(B(z, t)\big) < 2. t^{\frac 12\kappa'} \ \text { for }  \ \delta^{1-\ve_0'}< t< 1
$$
hence
$$
\operatornamewithlimits\supp\limits_z \mu_2 \big(B(z, t)\big) < 2. t^{\frac 14\kappa'} \ \text { for }  \ \delta< t< 1.\tag 8.9
$$

Regarding (7.12), we take $u_1=1, u_2=\frac{z_1}{z_2}$.
From (8.4)
$$
\max(|\Im z|, |\Im \frac {z_1}{z_2}z|)\gtrsim \rho  \ \text { if }  \ |z|=1
$$
which gives (7.2).

Applying Proposition 1$'$, either (8.5) or one of the following \hfill\break
(8.10) there is some $b\in\supp\mu$ such that
$$
N\Big(A +\frac {b-b'}{b''-b'}A, \delta\Big) > \delta^{-\sigma-\tau}.
$$
Then
$$
N\big((b'' -b')A+(b-b')A, \delta\big)> |b'-b''|^d  \ \delta^{-\sigma-\tau}> c\delta^{-\sigma -\tau+d\ve_0'}>\delta^{-\sigma-\frac\tau 2}
$$
for $\ve_0'$ small enough.

From the sumset inequalities
$$
N(A+bA, \delta)>\delta^{-\sigma -\frac\tau 8} \text { for some $b\in\supp \mu$}.
$$

\itemitem {(8.10)} $N(A+\frac{z_2}{z_1} A, \delta)> \delta^{-\sigma-\tau}$, implying
$$
N(z_1A+z_2A, \delta)> \delta^{d\ve_0'-\sigma-\tau}> \delta^{-\sigma-\frac\tau 2}
$$
and
$$
N(A+z_1A, \delta)+ N(A+z_2 A, \delta)>\delta^{-\sigma-\frac \tau 4}.
$$
This proves Proposition 2.

Iteration of Proposition 2 gives.

\proclaim
{Corollary 3}

Given $\sigma, \kappa, \kappa', \rho, \ve_1>0$, there are $\ve_0, \ve_0'>0$ and some $r\in\Bbb Z_+$ such that the following holds.

Let $\delta>0$ be small enough. Let $A, B\subset\Bbb C\cap B(0, 1)$ satisfy
\roster
\item "{\rm(8.12)}" $N(A, \delta)= \delta^{-\sigma}$.

\item "{\rm (8.13)}" $N\big(A\cap B(z, t),\delta\big) < t^\kappa N(A, \delta)$ for $\delta< t<\delta^{\ve_0}$ and $z\in\Bbb C$.

\item "{\rm (8.14)}" There is a probability measure $\mu$ on $B$ such that
$$
\mu\big(B(z, t)\big) < t^{\kappa'} \ \text { if } \ \delta<t<\delta^{\ve_0'} \text { and } z\in\Bbb C.
$$

\item "{\rm (8.15)}" There are elements $b_0, b_1, b_2 \in B$ such that
$$
|b_0- b_1|\sim|b_0-b_2|\sim \delta^{\ve_0'} \ \text { and } \ \Big|\Im \frac{b_0-b_1}{b_0-b_2}\Big|>\rho.
$$
\endroster

Then there are elements $z_1, \ldots, z_r$ obtained as product of at most $r$ elements from $B$, such that
$$
N(z_1A+\cdots+z_rA, \delta)> \delta^{-2+\ve_1}.\tag 8.16
$$
\endproclaim

Recall the following result from [B] (see Theorem 6 and its proof).

\proclaim
{Proposition 4} Let $\mu$ be a probability measure on $[0, 1]$ satisfying for some constants $\kappa>0, C$
$$
\mu(I) <C\rho^\kappa \ \text { if $I$ is a $\rho$-interval, $\delta<\rho<1$}\tag 8.17
$$
($\delta$ assumed small enough).

Then for some $s=s(\kappa, C)\in\Bbb Z_+$, the set $sA^{(s)}-sA^{(s)}$ where $A=\supp\mu$ and $sA^{(s)}= s$-fold sumset of $s$-fold product set $A^{(s)}$ of
$A$ is $\delta$-dense in $[0, 1]$.
\endproclaim

Note that in the conclusion of Proposition 4, we may clearly replace $\delta$ by any given power of $\delta$ (as a consequence of the statement)

\proclaim
{Proposition 5} 
Given $\kappa, \ve_1>0$, there is $\ve_0>0$ and $s\in\Bbb Z_+$ such that the following holds for $\delta>0$ small enough.

Let $A\subset\Bbb C\cap B(0, 1)$ satisfy
$$
N\big(A\cap B(z, t), \delta\big)< t^\kappa N(A, \delta)\tag 8.18
$$
for $z\in\Bbb C$ and $\delta<t<\delta^{\ve_0}$.

Then there is a time segment $T\subset \Bbb C$ of size at least $\delta^{\ve_1}$, such that each point in $T$ is $\delta$-close to an element from
$sA^{(s)}-sA^{(s)}$.
\endproclaim

Again in the conclusion, $\delta$ may be replaced by any fixed power of $\delta$.

\noindent
{\bf Proof.}

We may clearly assume $1\in A$, replacing $A$ by $\frac 1zA$ with $z\in A$ the element of largest norm.
Denote $\tilde A$ sets of elements obtained from $A$ by (boundedly many) sums of products.

Using Proposition 4, it is easily seen that it suffices to obtain a segment $T\subset\Bbb C$ of at least $\delta^{\ve_2}$ such that
$$
N(\tilde A\cap T_\delta, \delta)> \delta^{-1+\ve_2}\tag 8.19
$$
where $T_\delta$ denotes a $\delta$-neighborhood of $T$ and $\ve_2=\ve_2(\ve_1)$.

Following the proof of Proposition 2, we start specifying some point $z_0\in A$ and $\delta^{\ve_0}< t_0<1$ such that for all $\delta< t<t_0$, $z\in\Bbb C$
$$
N\big(A\cap B(z, t), \delta)< \Big(\frac t{t_0}\Big)^{\kappa/2} N\big(A\cap B(z_0, t_0),\delta\big).\tag 8.20
$$
By translation, we may assume $z_0=0$.
Performing another rescaling, we obtain a $\delta$-separated set $A\subset B(0, 1)$ such that
\roster
\item "{(8.21)}"  $0\in A, 1\in A$

\item "{(8.22)}" $|A\cap B(z, t)|< t^{\kappa/2}|A| \ \text { if } \ \delta< t< c, z\in\Bbb C$ .
\endroster
Define next
$$
|\eta|= \max_{z\in A} |\Im z|\tag 8.23
$$
and let
$$
z_1 =t_1+i\eta\in A.\tag 8.24
$$
If $|\eta|>c$, we may apply Corollary 3 with $A=B$ since (8.15) holds with $\ve_0'= 0, \rho =0(1)$ (take $b_0-0, b_1=1, b_2=z_1$).
From (8.16), $N(\tilde A, \delta)> \delta^{-2+\ve_2}$, implying (19) for some segment $T\subset\Bbb C$.

Assume $\eta=o(1)$.

Let $\nu$ be the image measure (on $[0, 1]$) of $\frac 1{|A|} 1_A$ under the map $z\to \Re z$.
It follows from (8.22), (8.23) that

\itemitem {(8.25)} $\nu(I)< \rho^{\kappa/3} \ \text { if } \ I\subset \Bbb R$ is a $\rho$-interval, $\eta<\rho<c$.

We apply Proposition 4 to $\nu$ (with $\delta$ replaced by $\eta$), concluding that $\widetilde{\supp\nu}$ and hence $\widetilde{\Re A}$ can be made $\eta$-dense in
$[0, 1]$.

Hence, by (8.23), this implies that
$$
[0, 1] \subset\tilde A+B(0, K\eta)\tag 8.26
$$
for some constant $K$.

In order to fulfill condition (8.15), we proceed as follows.

Take $z_2, z_3 \in\tilde A$ such that
$$
|\Im z_2|, |\Im z_3|< K\eta\tag 8.27
$$
$$
|\Re z_2- 10 Kt_1|< K\eta\tag 8.28
$$
$$
10K\eta< |\Re z_2- \Re z_3|< 11K\eta\tag 8.29
$$
\big (which is possible by (8.26)\big)
\input fig2.tex

If we let $b_0=2Kz_1, b_1=z_2, b_2=z_3$, then (8.15) clearly holds with $\rho =0(1)$ and
$$
\eta\equiv \delta^{\ve_0'}\tag 8.30
$$
for some $\ve_0'>0$.

We distinguish 2 cases.

If $\ve_0'$ in (8.30) is small enough, we may again apply Corollary 3 and conclude as in the case $\eta =0(1)$.

Otherwise, we proceed as follows.

Denote
$$
z_4 =10K z_1 -z_2\tag 8.31
$$
$$
A_1 =\{ z\in \tilde A; 0<\Re z< 1, |\Im z|< K\eta\}\tag 8.32
$$
$$
A_2 =z_4A_1.\tag 8.33
$$
By (8.26)
$$
[0, 1] \subset A_1+ B(0, K\eta)\tag 8.34
$$
$$
A_2\subset z_4 [0, 1]+B(0, 20 K^2 \eta^2)\tag 8.35
$$
and
$$
z_4[0, 1] \subset A_2+B(0, 20K^2\eta^2).\tag 8.36
$$
One easily verifies that
$$
\align
R=\{z=x+iy; 0\leq x\leq 1, 0\leq y\leq \eta\} &\subset A_1+ A_2 +\{z= x+iy; 0\leq x\leq c\eta, 0\leq y\leq c\eta^2\}
\\
&\subset \tilde A+ C\eta R.\tag 8.37
\endalign
$$
In particular, there is some $z_5\in \tilde A$ such that $\Re z_5 \sim \eta$ and
$\Im z_5 \sim\eta^2$.
From (8.37), also
$$
R\subset \tilde A+z_5 R\tag 8.38
$$
and hence, multiplying both sides of (8.38) by $z_5$
$$
\align
z_5 R& \subset \tilde A+z_5^2 R\\
R&\subset \tilde A+z_5^2 R.\tag 8.39
\endalign
$$
After a few iterations, we conclude that
$$
[0, 1] \subset\tilde A+B(0, \delta)\tag 8.40
$$
and in particular (8.19).

This completes the proof of Proposition 5.\hfill$\square$

\proclaim
{Proposition 6}

Given $\sigma>0$, there is $C(\sigma)>0$ and $s<s(\sigma)\in \Bbb Z_+$ such that for $\delta>0$ small enough the following holds.

Let $A\subset\Bbb C\cap B(0, 1)$ and
$$
N(A, \delta)> \delta^{-\sigma}.\tag 8.41
$$
Then there is a line segment $T\subset\Bbb C$ of length $\delta^\gamma$, 
$0<\gamma<C(\sigma)$ such that each point in $T$ is $\delta^{\gamma+\frac 12}$-close
to an element of $sA^{(s)} -sA^{(s)}$.
\endproclaim

\noindent
{\bf Proof.} Assume $A$ consists of $\delta$-separated points.

Take $t_0>\delta$ minimum such that for some $z_0\in A$
$$
|A\cap B(z_0, t_0)|\geq t_0^{\sigma/2}|A|.\tag 8.42
$$
Since obviously $|A\cap B(z_0, t_0)| \lesssim \big(\frac {t_0}\delta\big)^2$, it follows that 
$$
t_0\gtrsim \delta^{\frac {2-\sigma}{2-\frac \sigma 2}}.\tag 8.43
$$
From definition of $t_0$
$$
|A\cap B(z, t)|\leq (2t)^{\sigma/2} |A| \ \text { if } \ \delta< t\leq t_0 \ \text  { and } \
z\in\Bbb C.\tag 8.44
$$
Also, there is $z_1 \in A$ such that $|z_0-z_1|=t_0$.
Define
$$
A_1=(z_1-z_0)^{-1} \big((A-z_0)\cap B(0, t_0)\big).\tag 8.45
$$
From (8.42), (8.44), if $\frac \delta {t_0} <t<1$ and $z\in\Bbb C$
$$
\align
|A_1\cap B(z, t)|&\leq |A\cap B(z_0+(z_1 -z_0)z, tt_0)|\\
& \leq (2tt_0)^{\sigma/2}|A|\\
&\leq (2t)^{\sigma/2} |A_1|.\tag 8.46
\endalign
$$
Since (8.43)
$$
|A_1\cap B(z, t)|< t^\kappa|A_1| \ \text { for } \  \delta< t< 1\tag 8.47
$$
and
$$
\kappa =\frac \sigma{ 4-\sigma}\tag 8.48
$$
(note that $A_1$ consists of $\delta$-separated points).

Apply Proposition 5 to $A_1$ with $\kappa =\frac\sigma{4-\sigma}, \ve_1=\frac 12(\ve_0 =0)$ to obtain a segment $T_1$ of size $\delta^{\frac 12}$ such that
$$
T_1\subset sA_1^{(s)}- sA_1^{(s)}+B(0, \delta)\tag 8.49
$$
with $s=s(\sigma)\in\Bbb Z_+$.
From (8.45) the rescaling of $T_1$ by a factor $t_0$ gives a segment $T$ of size $\delta^{\frac 12} t_0^s=\delta^\gamma$ for which by (8.49)
$$
T\subset sA^{(s)} -sA^{(s)} +B(0, t_0^s\delta).\tag 8.50
$$
This proves Proposition 6.\qed

\noindent
{\bf Remark.}

From the statement of Proposition 6, it also follows that for any given integer $r\in\Bbb Z_+$, assuming (8.41), there is a segment $T\subset\Bbb C$ of length $\delta^\gamma$ such that
each point of $T$ is $\delta^{\gamma+r}$-close to an element from $sA^{(s)} -sA^{(s)}$, where $
\gamma=\gamma(\sigma, r)$ and $s=s(\sigma, r)$.

There is the following Cartesian version of Proposition 6 for $\Bbb C^d$ equipped with its product ring structure.
This is the result we need for our $SU(d)$-analysis (see Corollary 10).

\proclaim
{Proposition 7} Let $A\subset\Bbb C^d\cap B(0, 1)$ satisfying
$$
N(A, \delta)> \delta^{-\sigma}\tag 8.51
$$
for some $0<\sigma\leq d$.

Then there is a unit vector $\xi \in\Bbb C^d$ such that
$$
[0, \delta^\gamma].\xi \subset sA^{(s)} -sA^{(s)} +B(0, \delta^{\gamma+1})\tag 8.52
$$
for some $0\leq \gamma< C(d, \sigma)$ and $s\in\Bbb Z_+$, $s< s(d, \sigma)$.
\endproclaim

\noindent
{\bf Proof.}
Proceed by induction on the dimension $d$.

From Proposition 6 and the Remark, statement holds for $d=1$.

Next, assume Proposition 7 up to dimension $d$.
Let $A\subset\Bbb C^{d+1}\cap B(0, 1)$ and
$$
N(A, \delta)> \delta^{-\sigma}\tag 8.53
$$
for some $\sigma>0$.

We will denote again by $\tilde A$ sets of the form $sA^{(s)} -sA^{(s)}$ for varying $s$.
If $I\subset \{1, \ldots, d+1\}$, $\pi_I$ stands for the coordinate restriction.

Rearranging the coordinates, we may assume that $B=\pi_{\{1, \ldots, d\}} (A)$ satisfies
$$
N(B, \delta)> \delta^{-\frac d{d+1} \sigma} > \delta^{-\frac\sigma 2}.\tag 8.54
$$
From the induction hypothesis, there is a unit vector $\xi\in\Bbb C^d =[e_1, \ldots, e_d]$ such that
$$
[0, \delta^\gamma].\xi \subset sB^{(s)}-sB^{(s)}+B(0, \delta^{\gamma+1})\tag 8.55
$$
for some $0\leq\gamma \leq \gamma(\sigma), s\in\Bbb Z_+, s< s(d, \sigma)$.
Hence we may introduce a function
$$
\vp:[0, \delta ^\gamma]\to \tilde A\tag 8.56
$$
satisfying
$$
|\pi_{1, \ldots, d]} \vp(x) -x\xi|< \delta^{1+\gamma} \text { for } 0\leq x\leq \delta^\gamma.\tag 8.57
$$
We will distinguish several cases.

\noindent
{\bf Case 1.} $N(\pi_{d+1}\big(\vp([0, \delta^\gamma])\big), \delta^{2(1+\gamma)})<\delta^{-\frac 12}$.

Then there are elements $x_1, x_2\in [0, \delta^\gamma], |x_1-x_2|\gtrsim \delta^{\gamma+\frac 12}$ such that
\hfill\break
$|\pi_{d+1}\big(\vp(x_1)-\vp(x_2)\big)|< \delta^{2(1+\gamma)}$.
Hence, from (8.57), for $0\leq x< \delta^\gamma$
$$
\align
\vp(x)[\vp(x_1)-\vp(x_2)]&=\big(x\pi_{\{1, \ldots, d\}} \big(\vp(x_1)-\vp(x_2)\big).\xi, 0\big)+ 0\big(\delta^{1+\gamma}
|\vp (x_1)-\vp(x_2)|+ \delta^{2(1+\gamma)}\big)\\
&=x.\delta^{\gamma_1} \xi'+0(\delta^{1+\gamma+\gamma_1})\tag 8.58
\endalign
$$
where $\delta^{\gamma_1}=|\pi_{[1, \ldots , d]}\big(\vp(x_1) -\vp(x_2)\big)|\sim |x_1-x_2|\gtrsim \delta^{\gamma+\frac 12} $ and
$\xi'=\delta^{-\gamma_1}(\xi, 0)$ a unit vector in $\Bbb C^{d+1}$.
Therefore
$$
[0, \delta^{\gamma+\gamma_1}]\xi' \subset\tilde A+B(0, \delta^{\gamma+\gamma_1+1}).\tag 8.59
$$

\noindent
{\bf Case 2.}
$$
N(\pi_{d+1}\big(\vp([0, \delta^\gamma])\big), \delta^{2(1+\gamma)})\geq \delta^{-\frac 12}.
$$
In particular, $S=\pi_{d+1} (\tilde A)$ satisfies
$$
N(S, \delta^{2(1+\gamma)})> \delta^{-\frac 12}
$$
and an application of the $d=1$ result gives a segment $J\subset\Bbb C e_{d+1}$ of size $\delta^{\gamma_1}$ such that
$$
J\subset \tilde S +B(0, \delta^{\gamma_1+1}).\tag 8.60
$$

\noindent
{\bf Case 2.1} 
Assume $\vp$ approximatively linear in the sense that for all $x, y, x+y \in [0, \delta^\gamma]$,
$$
|\vp(x+y)-\vp(x)-\vp(y)|< \delta^{\gamma+\frac 12}.\tag 8.61
$$
Take $m=[\delta^{-\frac 14}]$ and $t_j = \frac jm \delta^\gamma(1\leq j\leq m)$.
Clearly (8.61) implies that
$$
\align
\vp(t_j)= &\vp(t_1)+\vp(t_{j-1})+0(\delta^{\gamma+\frac 12})\\
&= j\vp(t_1)+0(\delta^{\gamma+\frac 14}).\tag 8.62
\endalign
$$
Therefore
$$
\{j\vp(t_1); 1\leq j\leq m\} \subset\tilde A+B(0, \delta^{\gamma+\frac 14}).\tag 8.63
$$
Note that $\delta^{\gamma+\frac 14} \approx |t_1| \leq |\vp(t_1)|\lesssim \delta^{\frac 14}$.
Let $\xi'= \frac{\vp(t_1)}{|\vp(t_1)|} \in \Bbb C^{d+1}$.

It follows from (8.63) that
$$
[0, \lambda]\xi' \subset \tilde A+B(0, \delta^{1/4}\lambda)\tag 8.64
$$
with $\lambda= |\vp(t_1)|\delta^{-1/4}$.

\noindent
{\bf Case 2.2.}
Assume there are $x, y\in [0, \delta^\gamma], x+y\in [0, \delta^\gamma]$ such that
$$
|\vp(x+y)-\vp(x)-\vp(y)|> \delta^{\gamma+\frac 12}.\tag 8.65
$$
Denoting $\zeta =\vp(x+y)- \vp(x) -\vp(y)\in\tilde A$, it follows from (8.57) that
$$
|\pi_{[1, \ldots, d]}\zeta|\lesssim \delta^{1+\gamma}\tag 8.66
$$
hence
$$
|\zeta_{d+1}|> \delta^{\gamma+\frac 12}.\tag 8.67
$$
Taking $r\in\Bbb Z_+$ (an integer to specify), write
$$
\align
\tilde A &\supset \zeta^r \tilde A\supset (0, \zeta^r_{d+1} \tilde S)+ 0(\delta^{r(1+\gamma)})\\
&\underset {(8.60)} \to \supset (0, \zeta^r_{d+1} J)+ 0(\delta^{r(1+\gamma)}+|\zeta_{d+1}|^r \delta^{\gamma_1+1})\\
&=(0, \zeta^r_{d+1} J)+0(|\zeta_{d+1}|^r \delta^{\gamma_1+1})\tag 8.68
\endalign
$$
using (8.66), (8.67) and taking $r$ large enough.
Again (8.68) provides a segment $[0, \delta^{\gamma_2}]\xi'$, $\delta^{\gamma_2} =|\zeta_{d+1}|^r \delta^{\gamma_1}$, contained in $\tilde A+B
(0, \delta^{\gamma_2+1})$.

In summary, we certainly obtain a unit vector $\xi' \in \Bbb C^{d+1}$ such that for some \hfill\break
$0\leq\gamma'\leq C(d, \sigma)$
$$
[0, \delta^{\gamma'}]\xi' \subset\tilde A+B(0, \delta^{\gamma'+\frac 14})\tag 8.69
$$
\big(as there is only a $\delta^{1/4}$ gain in (8.64)\big).
Since the same statement holds with $\delta$ replaced by $\delta^4$ (note that then $\sigma$ in (8.41) needs to be replaced by $\frac\sigma 4$),
we proved (8.52) in $\Bbb C^{d+1}$.\qed

\proclaim
{Corollary 8}

Denote $\Bbb C^d \simeq\Delta\subset \Mat_{d\times d}(\Bbb C)$ the diagonal matrices.

Assume $A\subset\Mat_{d\times d}(\Bbb C)$ satisfies
\roster
\item"{(8.70)}" $A\subset B(0, 1)$

\item"{(8.71)}" $N(A, \delta)> \delta^{-\sigma}$

\item"{(8.72)}" $\dist (x, \Delta)< \delta$ for $x\in A$

Then there is $\xi \in \Delta, \Vert\xi\Vert =1$ such that

\item"{(8.73)}" $[0, \delta^\alpha]\xi\subset s' A^{(s)}- s'A^{(s)}+B(0, \delta^{\alpha+\beta})$
\endroster

\noindent
where $0\leq\alpha< c(d, \sigma), \beta>c(d, \sigma)>0$ and $s, s'\in\Bbb Z_+$; $s, s'< s (d, \sigma)$.
\endproclaim

\noindent
{\bf Proof.}

From (8.71), also $N(A, \delta^{1-\frac \sigma{2d}})> \delta^{-\frac \sigma 2}$.
Let $k\in\Bbb Z_+$ to be specified and
$$
\delta_1=\delta^{\frac 1k(1-\frac\sigma{2d})}.
$$
Clearly
$$
N(A, \delta^{1-\frac\sigma{2d}})\leq \prod_{0\leq\ell< k}\max_y N\big(A\cap B(y, \delta_1^\ell), \delta_1^{\ell+1}\big
).
$$
It follows that there is a subset $A_1\subset A$ and $\delta_1^{k-1} \leq \delta_2\leq 1$ such that

\roster
\item "{(8.74)}" $\diam A_1\leq \delta_2$

\item "{(8.75)}" $N(A_1, \delta_1\delta_2)> \delta^{-\frac \sigma{2k}}> \delta_1^{-\frac\sigma2}$.
\endroster

By (8.72), for each $x\in A$, there is $x'\in\Delta$ with $|x-x'|<\delta$; the set $B_1=\{x'; x\in A_1\}$ still satisfies (8.74), (8.75).

Take $\zeta \in B_1$ and denote
$$
B=\frac 1{\delta_2} (B_1-\zeta).\tag 8.76
$$
Hence
\roster
\item "{(8.77)}"  $B\subset \Delta\cap B(0, 1)$
\item "{(8.78)}"  $N(B, \delta_1)> \delta_1^{-\frac \sigma 2}$
\item "{(8.79)}" If $x\in B$, then $\dist (x, \frac{A_1-A_1}{\delta_2})\leq 2\frac\delta{\delta_2}.$
\endroster

Apply Proposition 7 to $B\subset\Bbb C^d$.
This gives a unit vector $\xi \in \Delta$ such that
$$
[0, \delta_1^\gamma]\xi \subset sB^{(s)} -sB^{(s)} +B(0, \delta_1^{1+\gamma})\tag 8.80
$$
for some $\gamma<\gamma (d, \sigma)$ and $s\in \Bbb Z_+, s< s(d, \sigma)$.

From (8.74), (8.79),
$$
\dist\Big(x, \frac{(A_1-A_1)^{(s)}}{\delta^s_2} \Big)< 2s \frac\delta{\delta_2} \ \text { for } \ x\in B^{(s)}
$$
and hence
$$
sB^{(s)} -sB^{(s)} \subset \delta_2^{-s} (s' A^{(s)}-s'A^{(s)})+B\Big(0, c_s\frac\delta{\delta_2}\Big).\tag 8.81
$$
From (8.80), (8.81)
$$
[0, \delta_1^\gamma \delta_2^s]\xi\subset s'A^{(s)} -s'A^{(s)}+B(0, \delta_1^{1+\gamma} \delta_2^s+c_s\delta\delta_2^{s-1}).\tag 8.82
$$
Note that by definition of $\delta_1$
$$
\frac\delta{\delta_2}< \delta^{\frac \sigma{2d}}<\delta_1^{1+\gamma}$$
provided
$$
k> 2d\big(1+\gamma(d, \sigma)\big)\sigma^{-1}.\tag 8.83
$$
Hence (8.73) holds with $\delta^\alpha =\delta_1^\gamma \delta_2^s$ and $\beta=\frac 1k\big(1-\frac\sigma{2d}\big)$.

\noindent
{\bf Remark.}

We also note that the element $\xi\in\Delta\simeq \Bbb C^d$ belongs to the algebra generated by $\{x'-y'; x, y \in A\}\subset\Bbb C^d$.

\bigskip

\centerline
{\bf Part II: Analysis on  the Unitary Group}

We now return to Theorem 1 and carry out the program sketched in the Introduction.

\bigskip
\centerline
{\bf \S9. Construction of Near-Diagonal Elements}

The main results from this section are formulated in Proposition 9 and Corollary 10.

Let $q\geq 2$ and $g_1, \ldots, g_q$ be fixed algebraic elements in $U(d)$ which freely generate the free group $F_q$. 
Let  $\Gamma=\langle g_1, \ldots,
g_q\rangle$. 

Denote by  $W_\ell$ the set of elements of $\Gamma$ that may be obtained as a word of length $\le \ell$ in $\{g_1, g_1^{-1}, \ldots, g_q, g_q^{-1}\}$.
Thus $W_\ell \subset W_{\ell+1}$. The following properties hold

\roster
\item"{(9.1)}" $|W_\ell|\sim (2q-1)^\ell$
\item "{(9.2)}" There is a constant $c_1>0$ such that $\Vert x-1\Vert>C^{-\ell}_1$ for $x\in W_\ell\backslash \{1\}$. (This is noncommutative diophantine property, introduced in [GJS] and used in connection with the  spectral gap theorem in [BG]).
\endroster

More generally, if $P(x_1, \ldots, x_r)$ is a polynomial with integer coefficients and variables $x_i\in \Mat_{d\times d}(\Bbb C)$, then

\roster
\item "{(9.3)}" If $x_1, \ldots, x_r\in W_\ell$, either $P(x_1, \ldots, x_r)=0$ or $|P(x_1, \ldots, x_r)|> C^{-\ell}$ where $C$ depends on $\Gamma$ and the degree
of $P$.
\endroster
\bigskip

For $\delta>0$, define
$$
W_{\ell, \delta} =\{x\in W_\ell:   \Vert x-1 \Vert<\delta\}.
$$
One may cover $U(d)$ by at most $\big(\frac {c}\delta\big)^{2d^2}$ balls $B_\alpha$ of size $\frac \delta 2$ and take $\alpha$ 
such that $|B_\alpha \cap W_{\frac\ell 2}|\gtrsim \delta^{2d^2}|W_{\frac \ell 2}|$.
Then $(B_\alpha\cap W_{\frac \ell 2})^{-1} (B_\alpha\cap W_{\frac\ell 2})\subset W_{\ell, \delta}$ and by (9.1)

\roster
\item"{(9.4)}" $|W_{\ell, \delta}|\gtrsim \delta^{2d^2} (2q-1)^{\frac\ell 2}$
\endroster

\medskip

The key idea underlying the proof of the following result originates in the  work of Helfgott [H].

\proclaim
{Lemma 9.5} Assume $A\subset U(d), b_1, \ldots, b_r \in W_\ell$ and $\delta>0$ with $\ell<\log \frac 1\delta$, such that
\roster
\item "{(9.6)}" $N(A, \delta)> \delta^{-\sigma}$

\item "{(9.7)}" span \, $A\subset$ \, span $(b_1, \ldots, b_r)$
\endroster

where ``span''  refers to the linear span in $\Mat_{d\times d}(\Bbb C)$.

Then there is $i \in \{1, \ldots, r\}$, $a\in A$ and a subset $A_1\subset (A\cup A^{-1} \cup\{b_1, b_1^{-1} , \ldots, b_r, b_r^{-1}\})^{(s)}$ (the $s$-fold
product set) and $\delta_1>0$ such that

\roster
\item "{(9.8)}" $\delta^C<\delta_1<\delta$ (where $C=C(\Gamma))$.

\item "{(9.9)}" $s\in\Bbb Z_+, s<s(\Gamma, \sigma)$.

\item "{(9.10)}" The elements of $A_1$ are $\delta_1$-separated and $|A_1|> \delta^{-\sigma/2r}$.

\item "{(9.11)}" $\Vert x ab_i^{-1} - ab_i^{-1} x\Vert<\delta_1$ for $x\in A_1$.
\endroster
\endproclaim

\noindent
{\bf Proof.}

We may assume that  $b_1, \ldots, b_r$  are linearly independent in $\Mat_{d\times d}(\Bbb C)$.

Consider the map
$$
\vp:\span (b_1, \ldots, b_r)\to \Bbb C^r: x\mapsto (\Tr x b_1^*, \ldots, \Tr x b_r^* ).
$$
Clearly, by (9.3)
$$
\align
\Vert\vp^{-1}\Vert& \sim \Vert[(\Tr b_ib_j^*)_{1\leq i, j\leq r}]^{-1}\Vert\\
&\lesssim |\det [(b_ib_j^*)_{1\leq i, j\leq r}]|^{-1}\\
&\lesssim C_2^{-\ell}.\tag 9.12
\endalign
$$
Let $A'\subset A$ be a $\delta$-separated set, $|A'|>\delta^{-\sigma}$.
By (9.12)
$$
|\vp(x) -\vp(y)|> C_2^{-\ell}\delta \ \text { for } \ x\not= y \ \text { in } A'
$$
and hence, for some $i=1, \ldots, r$
$$
N(\Tr A'b_i^{-1}, C_2^{-\ell}\delta)> \delta^{-\frac\sigma r}.\tag 9.13
$$
Take $A''\subset A'$ such that
\roster
\item "{\rm (9.14)}" $|A''|>\delta^{-\frac \sigma r}$

\item "{\rm (9.15)}" $|\Tr ab_i^{-1} -\Tr a' b_i^{-1}| >C_2^{-\ell}\delta$ for  $ a\not= a$ in $A''$.
\endroster

Denote by 
$$
\Cal A_s =(A\cup A^{-1}\cup \{ b_1, b_1^{-1}, \ldots, b_r, b_r^{-1}\})^{(s)}
$$
the $s$-fold product set and let  
$$
\delta_1 =C_2^{-\ell}\delta.\tag 9.16
$$
Since $\ell <\log \frac 1\delta$, trivially $N(\Cal A_s, \delta_1)\lesssim \big(\frac 1{\delta_1}\big)^{2d^2}< \big(\frac 1\delta\big)^{2d^2(1+\log C_2)}$. 
Hence, there  is some $s\in\Bbb Z_+$,
$$
\log s<4r d^2\sigma^{-1} (1+\log C_2)\tag 9.17
$$
such that
$$
N(\Cal A_{2(s+1)}, \delta_1)< \delta^{-\frac\sigma{2r}}N(\Cal A_s, \delta_1).\tag 9.18
$$
For $a\in A''$, consider the restricted conjugacy classes
$$
C_a =\{xab_i^{-1} x^{-1}, x\in\Cal A_s\}\subset\Cal A_{2s+2}.
$$
By (9.15), $\dist (C_a, C_{a'})> \delta_1$ for $a\not= a'$ in $A''$ and hence there is $a\in A''$ such that
$$
N(C_a, \delta_1) <N(\Cal A_{2s+2}, \delta_1) \delta^{\frac \sigma r}< N(\Cal A_s, \delta_1) \delta^{\frac \sigma {2r}}\tag 9.19
$$
\big(invoking (9.14), (9.18)\big).

Consider the map $\Cal A_s \to C_a: x\mapsto zab_i^{-1} x^{-1}$.
It follows from (9.19) that there is some $x_0\in\Cal A_s$ and a subset $\Cal A\subset\Cal A_s$ satisfying the following properties

\roster
\item "{\rm (9.20)}" $|\Cal A|>\delta^{-\frac\sigma{2r}}$ with $\delta_1$-separated elements.

\item "{\rm (9.21)}" $\Vert xab_i^{-1} x-x_0ab_i^{-1} x_0\Vert<\delta_1$ for $x\in\Cal A$.
\endroster

Hence the set $A_1=x_0^{-1} \Cal A$ satisfies 9.11.
This proves the Lemma.\qed

Our aim is to prove the following

\proclaim{ Claim $(*)$}\
Given $\delta>0$, there is $0<\delta_1<\delta$ and a subset $A\subset W_\ell$ such that

\roster
\item"{\rm(9.22)}" $\ell \sim \log \frac 1\delta\sim \log\frac 1{\delta_1}$
\item "{\rm(9.23)}" The elements of $A$ are $\delta_1$-separated and $|A|>\delta_1^{-c}$.
\item"{\rm(9.24)}" In an appropriate orthonormal basis
$$
\dist (x, \Delta)<\delta_1,
$$
where, as before,  $\Bbb C^d  \simeq  \Delta  \subset \Mat_{d\times d}(\Bbb C)$.
\endroster
\endproclaim

More generally, the same statement (with essentially the same proof) holds if $W_\ell$ is replaced by a large subset $H\subset W_\ell$,
i.e. $\log |H|\sim \ell$ (see the discussion preceding Proposition 9).

Compared with [BG2], it should be noted that the construction of the almost diagonal set $A$ does not use regular elements and this makes the
argument a bit more complicated.
On the other hand, Proposition 9 below gives a more general result of some  independent interest.

Assume $(*)$ fails for some (sufficiently small) $\delta$.
By induction on $1\leq d_1\leq d$, we then establish the following statement.

\proclaim
{$(**)$} If $\log \frac 1{\delta_1} \sim\log\frac 1\delta$, there is an element $x\in W_{\ell, \delta_1}$ with $\ell<C\log \frac 1\delta$ such
that $x$ has at least $d_1$ distinct eigenvalues.
\endproclaim

Next, we show that $(**)$ for $d_1=d$ implies $(*)$, hence obtaining a contradiction.

\noindent
{\bf Proof of $(**)\Rightarrow (*)$}

Take $x\in W_{\ell, \delta}, \log\frac 1\delta<\ell< C\log\frac 1\delta$ with $d$ distinct eigenvalues $\lambda_1, \ldots, \lambda_d$.
Since by (9.3)
$$
\prod_{i\not=j} |\lambda_i -\lambda_j|^2 \sim |\Res (P_x, P_x')|>C^{-\ell}
$$
where $P_x$ denotes the characteristic polynomial of $x$, it follows that
$$
|\lambda_i-\lambda_j| >C_3^{-\ell} \ \text { for } \ 1\leq i \not=  j\leq d.\tag 9.25
$$
Take $L>\ell, L\sim\ell$ and $L/\ell$ sufficiently large (according to the argument that follows).

Assume
$$
\span(W_{L, 2\delta})=\span(W_{\frac 12L, \delta_1})\tag 9.26
$$
where
$$
\delta_1 =(2C_3)^{-\ell}<\delta.\tag 9.27
$$
Take $b_1', \ldots, b_r'\in W_{\frac 12L, \delta_1}$, such that (9.26) = $\span (b_1', \ldots, b_r')$.
Note that $b_i=b_i'x\in W_{\frac 12L+\ell, \delta_1+\delta}\subset W_{L, 2\delta}$ and $b_1, \ldots, b_r$ are linearly independent.
Hence, by (9.26)
$$\span(b_1, \ldots, b_r) =\span (W_{L, 2\delta}).\tag 9.28
$$
Apply Lemma 9.5 with $A=W_{\frac 12L, \delta_1}$ and $\delta=C_1^{-L}$ \big(cf. (9.2)\big).
By (9.4)
$$
|W_{\frac 12 L, \delta_1}|\gtrsim \delta_1^{2d^2} (2q-1)^{\frac L4}> 3^{\frac L8}
$$
if $L>8d^2(\log C_3)\ell$.
Hence we may take $\sigma> \frac 1{10 \log C_1}$.

From the Lemma, we obtain $i=1, \ldots, r, a\in W_{\frac 12 L, \delta_1}$ and $A_1\subset W_{sL, 2s\delta}\subset W_{sL}$ such that for some
$C_1^{-L}>\delta_2> C_1^{-cL}$.

\roster
\item"{(9.29)}"  The elements of $A_1$ are $\delta_2$-separated and $|A_1|> 3^{\frac L{16r}}$.

\item "{(9.30)}" $\Vert yab_i^{-1}- ab_i^{-1} y\Vert < \delta_2$ for $y\in A_1$.
\endroster 

Note that $\xi =ab_i^{-1} =ax^{-1} (b_i')^{-1} =x^{-1}+O(\delta_1)$ so that the eigenvalues of $\xi$ satisfy essentially the same separation property
(9.25).
Choosing an orthonormal basis that makes $\xi$ diagonal, it follows from (9.30) and (9.25) that the off-diagonal elements of $y\in A_1$ are bounded
by $\delta_2C_3^\ell=\delta_3$.
Take a subset $A\subset A_1$ of $\delta_3$-separated elements, $|A|>C_3^{-2d^2\ell}|A_1|> 3^{\frac L{20r}} >\delta_3^{-c}$ 
\big(assuming $L>100 rd^2(\log C_3)\ell\big)$.

Hence $A$ satisfies $(*)$ (with $\ell =sL$ and $\delta_1=\delta_3$).

Next, assume that 9.26 fails, thus
$$
\dim\span (W_{L, 2\delta})> \dim\span(W_{\frac 12L, \delta_1}).
$$
Replace $L$ by $\frac L2$ and $\delta$ by $\frac{\delta_1}2$ and repeat the argument.

After at most $d^2$ steps (and assuming $L/\ell$ large enough), we reach
the same conclusion.
This proves the implication $(**)\Rightarrow (*)$.\qed
\medskip

\noindent
{\bf Proof of $(**)$ (assuming $(*)$ fails)}

For $d_1=1$, the statement is trivial.

For $d_1=2$, we may argue as follows.

If $T$ is a commutative subgroup of $GL_d(\Bbb C)$, then
$$
|W_\ell \cap T|\lesssim \ell \tag 9.31
$$
(using the fact that commutative subgroups of $F_q$, the free group on $q$ elements,  are cyclic).

Hence, if $\ell> c\log \frac 1{\delta_1}$, it follows from (9.4), (9.31) that
$$
W_{\ell, \delta_1} \not\subset \{z 1\!\!1; z\in\Bbb C, |z|=1\}
$$
and hence $W_{\ell, \delta_1}$ contains an element with at least $2$ distinct eigenvalues.

We now turn to the  inductive step.   Assume $(**)$ holds for $2\leq d_1< d$.
We follow the proof $(**)\Rightarrow (*)$.
Take $x\in W_{\ell, \delta}, \ell< c\log\frac 1\delta$ with $d_1$ distinct eigenvalues $\lambda_1, \ldots, \lambda_{d_1}$, hence satisfying

\roster
\item"{(9.32)}" $|\lambda_i-\lambda_j|> C_3^{-\ell} \ \text { for } \ 1\leq i\not= j\leq d_1$.
\endroster

Repeating the reasoning following $(9.26) $ (with the same notation), we obtain some $\xi\in W_{L+\ell, 2\delta_1+\delta}, \Vert\xi-x^{-1}\Vert
=O(\delta_1)$ and $A_1\subset W_{sL, 2s\delta}$ such that
\roster
\item "{(9.33)}" The elements of $A_1$ are $\delta_2$-separated and $|A_1|>\delta_2^{-c}$.

\item "{(9.34)}" $\Vert y\xi-\xi y\Vert  < \delta_2 $ for $y\in A_1$.
\endroster
where $\log \frac 1{\delta_2}\sim L$ and $C\ell <L<C'\ell$.

Take a subset $A_2\subset A_1A_1^{-1}\subset W_{2sL}$ satisfying
\roster
\item"{(9.35)}" $A_2 \subset B( 1\!\!1, C_3^{-2\ell})$.

\item "{(9.36)}" The elements of $A_2$ are $\delta_3 =C_3^{2\ell}. \delta_2$-separated and
$$
|A_2|> C_3^{-8d^2\ell}|A_1|>|A_1|^{\frac 12} >\delta_3^{-c}.
$$
\endroster
Obviously from (9.34)  we have 
\roster
\item"{(9.37)}" $\Vert y\xi-\xi y\Vert< 2\delta_2$ for $y\in A_2$.
\endroster

If $\xi$ has at least $d_1+1$ distinct eigenvalues, we are done.

Hence, assume $\xi$ has only $d_1$ distinct eigenvalues $\lambda_1', \ldots, \lambda_{d_1}'$.
Since $\Vert\xi-x^{-1}\Vert\lesssim \delta_1= (2C_3)^{-\ell}$, it follows from (9.32) that $\{\lambda_1', \ldots, \lambda'_{d_1}\}$ is an
$O(\delta_1)$-perturbation of $\{\lambda_1^{-1}, \ldots, \lambda_{d_1}^{-1}\}$ and in particular
$$
|\lambda_i'-\lambda_j'|> \frac 12 C_3^{-\ell} \ \text{ for } \ 1\leq i\not=j\leq d_1.\tag 9.38
$$
Diagonalize $\xi$ in a basis $e_1', \ldots, e_n'$ and write $\{1, \ldots, n\}=\bigcup_{s=1}^{d_1}I_s$ where
$$
\xi e_i'=\lambda_s' e_i' \ \text { for } \ i\in I_s.\tag 9.39
$$
We denote by $R_I$ the restriction to $I\subset \{1, \ldots, n\}$.

It follows from (9.38), (9.37) that
$$
\Vert y-y'\Vert\lesssim C_3^\ell \delta_2\tag 9.40
$$
where
$$
y'= \operatornamewithlimits\oplus\limits^{d_1}_{s=1} R_{I_s} y R_{I_s} = \operatornamewithlimits\oplus\limits^{d_1}_{s=1} y_s'.\tag 9.41
$$
For $s=1, \ldots, d_1$, let $y_s' =U_sP_s, P_s=\big((y_s')^*y_s'\big)^{\frac 12}$, be the polar decomposition of $y_s'$.
Since by (9.40)
$$
\Vert(y_s')^*y_s'- 1_{I_s}\Vert \lesssim C_3^\ell \delta_2\tag 9.42, 
$$
it follows that
$$
\Vert y_s' -U_s\Vert\lesssim C_3^\ell \delta_2.\tag 9.43
$$
We distinguish 2 cases.

\roster
\item "{(i)}" Assume that for all $y\in A_2$ and $s=1, \ldots, d_1$
$$
\dist(y_s', \{z 1\!\!1_{I_s}; z\in\Bbb C\})<\delta_3.\tag 9.44
$$
From (9.40), $y$ is an $O(\delta_3)$-perturbation of a diagonal matrix.

Hence $A_2$ satisfies $(*)$ (with $\ell=2sL, \delta_1=\delta_3)$, a contradiction.
\medskip

\item "{(ii)}" Let $y\in A_2$ and $s=1, \ldots, d_1$ (say $s=1$) such that
$$
\dist (y_1' , \{z 1_{I_1}; z\in\Bbb C\})\geq \delta_3.\tag 9.45
$$
Hence, by (9.43)
$$
\dist(U_1, \{z 1_{I_1}; z\in\Bbb C\})>\frac 12\delta_3\tag 9.46
$$
and $U_1$ has at least 2 eigenvalues that are $\frac 12 \delta_3$-separated.
\endroster

Take for $s=1, \ldots, d_1$ a basis $\{e_i''; i\in I_s\}$ of $[e_i'; i\in I_s]$ diagonalizing $U_s$.
If
$$
U_s e_i'' =\mu_i e_i'' \text { for } \ i\in I_s\tag 9.47
$$
we have
$$
\align
|1-\mu_i|\leq \Vert 1_{I_s} -U_s\Vert&\underset{(9.43)}\to < \Vert 1_{I_s} - y_s'\Vert+ 0(C_3^\ell \delta_2)\\
&< \vert 1-y\Vert+\Vert y-y'\Vert+0(C_3^\ell\delta_2)\\
& \operatornamewithlimits < \limits_{\Sb (9.35)\\ (9.40)\endSb}
C_3^{-2\ell}+0(C_3^\ell \delta_2)< 2C_3^{-2\ell}\tag 9.48
\endalign
$$
and
$$
\big(\xi.(\oplus_s U_s)\big)(e_i'')=\mu_i\lambda_s' e_i'' \ \text { for } \ i\in I_s.\tag 9.49
$$
By assumption, we may take $\mu_1, \mu_2, \{ 1, 2\}\subset I_1$, such that
$$
|\mu_1-\mu_2|>\frac 12\delta_3.\tag 9.50
$$
From (9.38), (9.48), we have for $i_1\in I_{s_1}, i_2\in I_{s_2}, s_1\not= s_2$
$$
|\mu_{i_1} \lambda_{s_1}' -\mu_{i_2} \lambda_{s_2}'|> \frac 12 C_3^{-\ell}- 4 C_3^{-2\ell}> \frac 13 C_3^{-\ell}.\tag 9.51
$$
In view of (9.50), (9.51), $\xi(\oplus_s U_s)$ has at least $d_1+1$ eigenvalues that are $\frac 12\delta_3$-apart.

Consider the element $\xi y\in W_{(2s+1)L+\ell, 3\delta_1+\delta}\subset W_{2(s+1)L, 2\delta}$.
Since \break $\Vert\xi y -\xi(\oplus U_s)\Vert~\lesssim~C_3^\ell\delta_2$ 
by (9.40), (9.43) and $\delta_3=C_3^{2\ell}. \delta_2$, also
$\xi y$ has at least $d_1+1$ distinct eigenvalues.
Hence $(**)$ holds for $d_1+1$.\qed

This completes the proof of $(*)$.
\bigskip

We also need the following extension.

Assume $H\subset W_\ell, 1\in  H =H^{-1}$ such that
$$
\log |H|\sim\ell.\tag 9.52
$$
The $s$-fold product set $H^{(s)}$ obviously satisfies $H^{(s)} \subset W_{s\ell}$.

On the other hand, from Razborov's product theorem in the free group (see [R])
$$
|H.H.H|\gg |H|^{2-\ve} \qquad (\ve>0).\tag 9.53
$$
Hence, for $\alpha<\frac{\log 2}{\log 3}$
$$
|H^{(s)}|> |H|^{s^\alpha}.\tag 9.54
$$
(Note that any statement of the form $\frac{\log |H^{(s)}|}{\log |H|} \longrightarrow \infty$ with $s\rightarrow \infty$ suffices
for our purpose.)

Replacing the sets $W_{\ell'}, \ell'<\ell$ by product sets $H^{(s)}$, a straightforward adaptation of previous analysis permits us to obtain
again a set $A\subset H^{(s)}$, for some $s$, satisfying (9.23), (9.24).
This gives

\proclaim
{Proposition 9}
Let $g_1, \ldots, g_q$ be algebraic elements in $U(d)(q\geq 2)$ generating a free group.
Take $H\subset W_\ell(g_1, \ldots, g_q)$, $\ell$ sufficiently large, such that
$$
\log|H|\sim\ell.\tag 9.55
$$
There is $A\subset H^{(s)}$ $(s<C)$ and $\delta>0$ such that
\roster
\item "{\rm (9.56)}" $\log \frac 1\delta\sim\ell$.

\item "{\rm (9.57)}" The elements of $A$ are $\delta$-separated and $|A|>\delta^{-c}$.

\item "{\rm (9.58)}" In an appropriate orthonormal basis,
$$
\dist(x, \Delta) <\delta \ \text { for } \ x\in A.\tag 9.59
$$
\endroster
\endproclaim

\proclaim
{Corollary 10} Let $g_1, \ldots, g_q \in SU(d) (q\geq 2)$ be algebraic and free.
Let $H\subset W_\ell (g_1, \ldots, g_q)$, $\ell$ large enough, such that
\roster
\item  "{\rm (9.60)}" $\log |H|\sim \ell$.
\endroster
Then there are $\delta_0> \delta>0$ with
\roster
\item"{\rm (9.61)}" $\log \frac 1{\delta_0} \sim \log \frac 1\delta\sim \ell$
\endroster
and $\xi=(\xi_{ij})_{1\leq i, j\leq d}$ with
\roster
\item"{\rm (9.62)}" $\xi_{ji} =\overline{\xi_{ij}}, \xi_{ii}= 0 $ and $\Vert\xi\Vert =1$
\endroster
such that the following holds.

Let  $\eta\in \Cal L(\Bbb C^d), \Vert\eta \Vert<\delta$ and $t\in[0, \delta_0]$. Then
$$
\Big\Vert 1+ t\sum_{i, j} \xi_{ij} \eta_{ij}(e_i\otimes e_j)-x\Big\Vert< \delta_0^{1+\gamma}\Vert\eta\Vert\tag 9.63
$$
for some $x\in (H\cup\{1+\eta, (1+\eta)^{-1}\})^{(s)}$, $s<C$ and where $\gamma>0$ is a fixed constant.
\endproclaim

\noindent
{\bf Proof.}

First apply Proposition 9 to $H$.
We obtain $A\subset H^{(s_1)}$, $s_1<C$ and $\delta_1>0$, $\ell\sim\log \frac 1{\delta_1}$ such that the elements of $A$ are $\delta_1$-separated
$$
\log|A|\sim\ell\tag 9.64
$$
and (after a base change)
$$
\dist (x, \Delta) <\delta_1 \ \text { for } x\in A.\tag 9.65
$$
Denote $V$ the vector space $M_{d\times d}(\Bbb C)$ for $x\in S U(d)$ and consider the adjoint representation $\rho_x, \rho_x(z)= x^{-1} zx$,
acting unitarily on $V$.

We will apply Corollary 8 to the set
$$
\Cal A=\{\rho_x; x\in A\}\subset U(V).\tag 9.66
$$
To each $x\in A$, associate
$$
x'=\sum x_{ii} \, e_i\otimes e_i\in\Delta
$$
for which by (9.65)
$$
\Vert x-x'\Vert \lesssim \delta_1.\tag 9.67
$$
Since $\det x=1$, it follows that
$$
\Big|1-\prod |x_{ii}|\Big|\leq \Big|1-\prod x_{ii}\Big|\lesssim \delta_1\tag 9.68
$$
and
$$
\big|1-|x_{ii}|\big|\lesssim \delta_1 \qquad (1\leq i\leq d).\tag 9.69
$$
Also
$$
\Vert\rho_x-\rho_{x'}\Vert\lesssim \delta_1\tag 9.70
$$
where $\rho_{x'} \in\Delta_V =  \{$diagonal elements of $\Mat(V)$\}.

For $x, y\in A$, we have
$$
\align
\Vert \rho_{x'} -\rho_{y'}\Vert &\sim \max_{i\not=j}\Big|\frac {x_{ii}}{x_{jj}}-\frac {y_{ii}}{y_{jj}}\Big|\\
&\gtrsim \max_i \Big|\frac {x_{ii}^{d-1}}{\prod_{j\not=i} x_{jj}} - \frac {y_{ii}^{d-1}}{\prod_{j\not= i} y_{jj}}\Big|\\
&\gtrsim \max_i|x^d_{ii} - y^d_{ii}|- O(\delta_1) \qquad (\text {by  (9.68)}).
\endalign
$$
Hence, if $\Vert\rho_{x'} -\rho_{y'} \Vert <\delta_1$, there are $k_i \in\{ 0, 1, \ldots, d-1\}$ $(1\leq i\leq d)$ such that
$$
\Big\Vert y-\sum_i e\Big(\frac {k_i}{d}\Big) x_{ii} (e_i\otimes e_i)\Big\Vert\lesssim\delta_1.
$$
Since the elements of $A$ are $\delta$-separated, we may find a subset $A_1\subset A$
such that $|A_1|\sim |A|$ and
$$
\Vert\rho_{x'} -\rho_{y'}\Vert> C\delta_1 \ \text { for } \ x, y\in A_1
$$
hence, by (9.70)
$$
\Vert\rho_x- \rho_y\Vert>\delta_1 \ \text { for } \ x, y\in A_1.
$$
It follows that
$$
N(\Cal A, \delta_1) \gtrsim |A|> \delta_1^{-c}.
$$
Thus $\Cal A$ satisfies (8.70)-(8.72) with $\delta$ replaced by $\delta_1$.

By Corollary 8, we obtain $\xi =(\xi_{ij})\in \Delta_V$, $\Vert\xi\Vert=1$ such that (8.73) holds.

Thus for $t\in [0, \delta_1^\alpha]$, there is $M\in s'\Cal A^{(s)} -s' \Cal A^{(s)}$ such that
$$
\Big\Vert t\sum \xi_{ij} z_{ij} (e_i\otimes e_j) -M(z)\Vert< \delta_1^{\alpha+\beta}\tag 9.71
$$
for all $z\in V$, $\Vert z\Vert \leq 1$.

Moreover (by the Remark following Corollary 8), $\xi$ belongs to the algebra generated by $\{\rho_{x'}-\rho_{y'}; x, y \in A\}$.
since $\rho_{x'}\in\Delta_V$ has diagonal elements $\frac{x_{ii}}{x_{jj}}$, it follows that $\xi_{ii}=0$.

Also 
$$
\langle \rho_x(e_i\otimes e_j), e_i\otimes e_j\rangle = \overline{\langle\rho_x(e_j\otimes e_i), e_j\otimes e_i\rangle}
$$
for $x\in U(d)$, hence
$$
\langle M(e_i\otimes e_j), e_i\otimes e_j\rangle =\overline{\langle M(e_j\otimes e_i), e_j\otimes e_j\rangle}.
$$
Therefore we may take $\xi$ with $\xi_{ji} =\overline{\xi_{ij}}$ in (9.71).

Note that $\Cal A^{(s)} =\rho_{A^{(s)}}$ and $\rho_{A^{(s)}} \eta =\{x^{-1} \eta x; x\in A^{(s)}\}$.

If $\Vert\eta\Vert <\delta$, it follows that
$$
\align
M\eta+1&=\prod^{s'}_{\alpha=1} (1+x_\alpha^{-1} \eta x_\alpha)\prod^{2s'}_{\alpha=s'+1} (1-x_\alpha^{-1} \eta x_\alpha)+O(\delta\Vert\eta\Vert)\\
&= \prod^{s'}_{\alpha=1} x_\alpha^{-1} (1+\eta)x_\alpha \prod^{2s'}_{\alpha=s'+1} x_\alpha^{-1} (1+\eta)^{-1} x_\alpha +O(\delta\Vert\eta\Vert)
\endalign
$$
for some $x_\alpha \in A^{(s)}$. Hence
$$
M\eta+1\in (H\cup\{1+\eta, (1+\eta)^{-1}\})^{2s'(2ss_1+1)}+B(0, C\delta\Vert\eta\Vert)
$$
and there is $x$ in a product set of $H\cup\{1+\eta, (1+\eta)^{-1}\}$ such that
$$
\Vert M\eta+1 -x\Vert\lesssim \delta\Vert\eta\Vert.\tag 9.72
$$
Taking $z =\frac\eta{\Vert\eta\Vert} $ in (9.71), it follows that
$$
\Vert 1+t\sum\xi_{ij} \eta_{ij}(e_i\otimes e_j)-x\Vert \lesssim \delta_1^{\alpha+\beta}\Vert\eta\Vert+\delta\Vert\eta\Vert.\tag 9.73
$$
Take $\delta<\delta_1^{\alpha+\beta}$.
We obtain (9.63) with $\delta_0=\delta_1^\alpha, \gamma=\frac\beta\alpha$.

This proves Corollary 10. \qed

\bigskip

\centerline
{\bf \S10. Expansion in $SU(d)$ \quad (1)}

Let $g_1, \ldots, g_k \in G=SU(d)$ be algebraic elements and assume $\Gamma =\langle g_1, \ldots, g_k\rangle$ Zariski dense in $G$. Denote
$$
\nu=\frac 1{2k}\sum^k_{i=1} (\delta_{g_i} +\delta_{{g_i}^{-1}})\tag 10.1
$$
which is a symmetric probability measure on $G$.

Our aim is to establish a spectral gap, i.e.
$$
\Vert f*\nu\Vert_2 \leq (1-\ve)\Vert f\Vert_2 \ \text { for } \ f\in L^2_0 (G).\tag 10.2
$$
Invoking a result of Breuillard and Gelander [BrGe], we may assume $k=2$ and $\{g_1, g_2\}$ generate the free group $F_2$.

Note that $G$ is $d_1=d^2-1$ dimensional over $\Bbb R$.

Denote for $\delta>0$
$$
P_\delta=\frac { 1_{B(1, \delta)}}{|B(1, \delta)|}\tag 10.3
$$
an approximate identity.

As we show below, (10.2) is a consequence of the following  main proposition.

\proclaim
{Proposition 11} For any given $\tau>0$, there is a positive integer $\ell<C(\tau)\log\frac 1\delta$ such that
$$
\Vert\nu^{(\ell)} *P_\delta\Vert_\infty < \delta^{-\tau}\tag 10.4
$$
where $\nu^{(\ell)} =\underbrace{\nu*\cdots*\nu}_{\ell}$ denotes the $\ell$-fold convolution.
\endproclaim

Note that (10.4) may be replaced by the a priori weaker statement
$$
\Vert\nu^{(\ell)}*P_\delta\Vert_2 < \delta^{-\tau}.\tag 10.5
$$
Indeed, since $P_\delta< C(P_\delta *P_\delta)$, one has for $x\in G$
$$
(\nu^{(2\ell)} *P_\delta)(x) \leq C \langle
\nu^{(\ell)} * P_\delta, \tau_{x^{-1}} (\nu^{(\ell)}*P_\delta)\rangle \leq\Vert \nu^{(\ell)}* P_\delta\Vert^2_2
$$
where $(\tau_{x^{-1}} f)(y)=f(yx^{-1})$.

Note that since $\{g_1, g_2\}$ are algebraic and free, one gets trivially that
$$
\Vert\nu^{(\ell)} *P_\delta\Vert_2 < \delta^{-\frac 12 d_1+\theta}\tag 10.6
$$
for some $\theta>0$.
Proposition 11 will therefore follow from a bounded number of applications of

\proclaim
{Lemma 10.7} ($L^2$-flattening lemma).

Given $\gamma>0$, there is $\kappa>0$ such that if $\delta>0$ is small enough, $\ell\sim \log\frac 1\delta$, if
$$
\Vert\nu^{(\ell)} *P_\delta\Vert_2  > \delta^{-\gamma}\tag 10.8
$$
then
$$
\Vert\nu^{(2\ell)}* P_\delta\Vert_2 <\delta^\kappa\Vert\nu^{(\ell)}* P_\delta\Vert_2.\tag 10.9
$$
\endproclaim

We use notation $| \ |$ to denote either the Haar measure on $G$ or the cardinality of a discrete set.

Denote $\mu=\nu^{(\ell)}*P_\delta$ and assume (10.9) fails, i.e.
$$
\Vert \mu*\mu\Vert_2 >\delta^{0+}\Vert\mu\Vert_2.\tag 10.10
$$
From a noncommutative version of the BSG theorem due to Tao [T, Thm. 5.4.], we obtain a subset $H$ of $G$, $H$ a union of $\delta$-balls, and a
discrete set $X\subset G$ satisfying the following properties.

\roster
\item"{(10.11)}" $H=H^{-1}$.
\item "{(10.12)}" $H.H\subset H.X\cap XH$.
\item "{(10.13)}" $|X|<\delta^{0-}$
\item "{(10.14)}" $\mu(aH)>\delta^{0+}$ for some $a\in G$.
\item "{(10.15)}" $|H|<\delta^\gamma$.
\endroster

Properties (10.11)- (10.13) mean that $H$ is an `approximate group'.
Note that (10.12), (10.13), (10.15) imply
\roster
\item "{(10.16)}" $|H^{(s)}|<_{(s)} \delta^{0-} |H|<\delta^{\gamma-} $ for any given $s\in\Bbb Z_+$
\endroster
where $H^{(s)}$ stands for the $s$-fold product set.

For notational simplicity, let $H'$ denote product sets $H^{(s)}$ for unspecified (but bounded) $s\in\Bbb Z_+$.

The proof of Lemma 7, and hence of  Proposition 11, will be completed by showing that there is no approximate group $H$ satisfying (10.14), (10.15).

This is the same approach as in [BG] for $SU(2)$.
But,  as discussed in the introduction, the argument used here differs from that of [BG],  and is analogous to [BG2] on expansion in $SL_d(p^n)$.

We first show how to derive (10.2) from Proposition 11.

In [BG] treating $G=SU(2)$, we relied on [GJS] extension of   Sarnak-Xue technique,  based on  suitable averaging  of characters of the irreducible representations on spaces of homogeneous
polynomials.
The argument presented below gives an alternative approach that is perhaps more geometric  and `general' in the sense of being applicable to other groups.

We will reduce the problem to a convolution property on the group $G$ (Lemma 10.35 below).
The relevant inequality will then be established first for $G=SU(2)$ and next in general for $G=SU(d)$, using $SU(2)$-embeddings.

Let $G=SU(d)$ and denote $(\rho_gf)(x) =f(xg)$ acting on $L^2(G)$.
 Letting $\nu$ be the discrete, symmetric probability measure (10.1), we have to
establish (10.2).
Assume to the contrary that
$$
\Vert f*\nu\Vert_2> 1-\ve\tag 10.18
$$
where $f\in L^2_0(G), \Vert f\Vert_2=1$.

We introduce  a Littlewood-Paley decomposition on $G$ (which is standard construction in harmonic analysis, see, for example, Chapter 4 in [S]).

For $f\in L^2(G)$, $0<\delta<1$, let
$$
(f*P_\delta)(x)= \nint_{B(1, \delta)} f(xg)dg \tag 10.19
$$
with $\Cal P_\delta$ introduced in (10.13) and denote
$$
\align
\Delta_1 f&= (f*P_{2^{-1}})-\int_Gf\\
\Delta_k f&= (f*P_{2^{-k}})-(f*P_{2^{-k+1}})\  \text { for } \ k\geq 2.\tag 10.20
\endalign
$$
Thus we have dcomposition
$$
f=\sum_{k\geq 1} \Delta_kf\tag 10.21, 
$$
satisfying the square function property
$$
\frac 1C \Big(\sum \Vert\Delta_k f\Vert_2^2\Big)^{\frac 12}\leq \Vert f\Vert_2\leq C\Big(\sum\Vert\Delta_k f\Vert^2_2\Big)^{\frac 12}.\tag 10.22
$$
Since $\rho_g(f*P_\delta)=(\rho_gf)*P_\delta$, it follows that $\Delta_k f*\nu=\Delta_k(f*\nu)$.
Hence, by (10.22), for all $\ell\in \Bbb Z_+$
$$
\Vert f*\nu^{(\ell)}\Vert_2 \sim\Big(\sum_{k\geq 1} \Vert (\Delta_kf) *\nu^{(\ell)}\Vert^2_2\Big)^{\frac 12}.\tag 10.23
$$
Note that from the Zariski-density assumption for $\Gamma$, no elements of $L_0^2(G)\backslash \{0\}$ are \break
$\Gamma$-invariant.
Hence
\bigskip

\roster
\item"{(10.24)} " $f\to 0$ (weakly) if $f\in L_0^2(G), \Vert f\Vert_2=1$ satisfies (10.18) with $\ve\to 0$.
\endroster

\medskip
Therefore we may assume
$$
f=\sum_{k\geq k_0} \Delta_kf\tag 10.25
$$
where $k_0=k_0(\ve)\overset{\ve\to 0}\to \longrightarrow \infty$.

Let $\ell_0\in \{2^m\}$ be fixed and sufficiently large (to be specified).

From (10.18), (10.22), (10.23) and (10.25), taking $\ve$ small enough
$$
\Vert f*\nu^{(\ell_0)}\Vert_2>(1-\ve)^{\ell_0}> \frac 12
$$
and
$$
\sum_{k>k_0} \Vert\Delta_k f*\nu^{(\ell_0)}\Vert_2^2> c
\sum_{k\geq k_0} \Vert\Delta_k f\Vert_2^2\tag 10.26
$$
with $c$ independent of $\ell_0$.

Therefore there is some $k>k_0$ such that
$$
F=\frac{\Delta_kf}{\Vert \Delta_k f\Vert_2}\tag 10.27
$$
satisfies
$$
\Vert F* \nu^{(\ell_0)}\Vert_2 > c.\tag 10.28
$$
Take $\ell\sim k$ (to specify). From (10.28)
$$
\Vert F*\nu^{(\ell_0\ell)} \Vert_2> c^\ell.\tag 10.29
$$
Let $\delta=4^{-k}$.
Then, by (10.27), $F\approx F*P_\delta$ and (10.29) implies
$$
\Vert F*\mu\Vert_2> c^\ell\tag 10.30
$$
where $\mu=\nu^{(\ell_0\ell)}*P_\delta$.

Fix $\tau>0$. If we take
$$
\ell_0\ell>C(\tau)\log\frac 1\delta
$$
(hence
$$
\ell> \frac{C(\tau)k}{\ell_0}\tag 10.31),
$$
Proposition 11 gives
$$
\Vert \mu\Vert_\infty <\delta^{-\tau}.\tag 10.32
$$
Rewrite (10.30) as 
$$
\int\limits_G \Big|\int\limits_GF(xg) \mu(g) dg\Big|^2 dx >c^{2\ell}
$$
and hence by (10.32)
$$
\delta^{-2\tau}\int\limits_G\Big|\int\limits_G F(xg)\overline{F(x)} dx\Big|dg \geq \iint\limits_{GG}\Big|\int\limits_G 
F(xh^{-1} g)\overline{F(x)} dx\Big| \mu(g)\mu(h) dgdh> c^{2\ell}.\tag 10.33
$$

Denoting $F_1(x) =\overline{F(x^{-1})}$, it follows from (10.33) that
$$
\Vert F_1* F\Vert_1 > c^{2\ell} \delta^{2\tau}> \delta^{\frac {C(\tau)}{\ell_0}+2\tau}.\tag 10.34
$$
In order to obtain a contradiction, it will therefore suffice to apply

\proclaim
{Lemma 10.35} Given $c>0$, there is $c'>0$ such that if $\delta_1>0$ is small enough and $F_1, F_2\in L^2 (G)$ satisfy $\Vert F_1\Vert_2 \leq 1, \Vert
F_2\Vert_2\leq 1$ and
$$
\Vert F_2 * P_{\delta_1}\Vert_2 <\delta_1^c.\tag 10.36
$$
Then
$$
\Vert F_1*F_2\Vert_1< \delta_1^{c'}.\tag 10.37
$$
\endproclaim

Indeed, by (10.27), $F$ satisfies (10.36) with $\delta_1= 2^{-k/2}$ say and (10.37) contradicts (10.34) by taking first $\tau$ small enough and then $\ell_0$
large enough.

Returning to Lemma 10.35 for $G=SU(d)$, note that this principle obviously fails for $d=1$.

We will first establish (10.35) for $G= SU(2)$ and then derive from this the statement for $G= SU(d), d>2$.

\noindent
{\bf Proof of Lemma 10.35 for $G=SU(2)$}

Denoting by $T$ the convolution operator by $F_2$ acting on $L^2(G)$, we have to prove that
$$
\Vert T\Vert <\delta_1^{c'}.\tag 10.38
$$
It suffices to verify (10.38) for the action of $T$ in the irreducible unitary representations of $SU(2)$, that is on the spaces $W_n =[z^k w^{n-k}; 0\leq k\leq
n]$ of homogeneous polynomials on the unit sphere of $\Bbb C^2$.
From the theory of induced representations (Frobenius' theorem) applied to the diagonal subgroup 
$D=\{R_\theta =\pmatrix e^{2\pi i \theta} & 0\\ 0&e^{-2\pi i\theta}
\endpmatrix ; 0\leq\theta <1\}$, we see that the IUR of $G$ are contained in one of the following two representations $\rho_0, \rho_1$ (depending on $n$ being 
even or odd).

(i) $\rho_0$ is acting by right translation on the subspace $V_0$ of $L^2(G)$ of functions that are left invariant under the action of $D$.
Thus $f\in V_0$ if $f(x) =f(R_\theta x)$ and $f$ factors over $D\backslash G \simeq S^{(2)}$.
Equivalently, $\rho_0$ may be seen as the representation of $SO(3)$ on $L^2(S^{(2)})$ by rotation.

(ii) $\rho_1$ is acting by right translation on the subspace $V_1$ of $L^2(G)$ of functions satisfying
$$
f(R_\theta x)= e^{i\theta} f(x).\tag 10.39
$$

To establish (10.38), we need to show that
$$
\Vert T|_{V_0} \Vert <\delta_1^{c'} \text { and } \Vert T|_{V_1} \Vert< \delta_1^{c'}.\tag 10.40
$$
We treat the action on $V_1$ (the case of $V_0$  is analogous).

Let $f\in V_1, \Vert f\Vert_2 =1$.
We have to prove that $\Vert f*F_2\Vert_2 <\delta_1^{c'}$.
In view of (10.36), we can assume that
$$
\Vert f*P_{\delta_2} \Vert_2 < \delta_2\tag 10.41
$$
with $\delta_2= \delta_1^{c'}$. Write
$$
\align
\Vert Tf\Vert_2^2 &=\iint\limits_{G\times G} \Big[\int\limits_G f(xy^{-1}_1) \overline{ f(xy_2^{-1})} dx\Big] F(y_1)\overline{F(y_2)} dy_1 dy_2\\
&\leq \Big\Vert\int\limits_G f(x) \overline{f(x.)} dx\Big\Vert_2.
\endalign
$$

Squaring again and using that $f\in V_1$, it follows
$$
\align
\Vert Tf\Vert^4_2&\leq \int\limits_{G\times G\times G} f(x) \overline{f(xy)} \ \overline{f(x_1)} f(x_1 y) dxdx_1dy\\
&=\int\limits_{G\times G\times G} \int^1_0 f(x) \overline {f(R_\theta xy)} \ \overline {f(x_1)} \ f(R_\theta x_1y) dxdx_1dyd\theta\\
&= \int\limits_{G\times G\times G} f(x) \overline {f(y)} \ \overline {f(x_1)} \Big[\int_0^1 f(R_\theta x_1 x^{-1} R_{-\theta}y)d\theta\Big]
\endalign
$$
and
$$
\align
\Vert Tf\Vert^8_2 &\leq \int_G \Big\Vert\int^1_0 f(R_\theta x R_{-\theta}.) d\theta \Big\Vert^2_2 \, dx\\
&=\int_{D\backslash G}\Vert S_x f\Vert^2_2 \, dx\tag 10.42
\endalign
$$
where we denoted
$$
S_zf(x) =\int^1_0 f(R_\theta z R_{-\theta} x)d\theta.\tag 10.43
$$
For $z =\pmatrix r& se^{i\psi}\\ -se^{-i\psi}& r\endpmatrix, s\not=0$, the operator $S_z$ acting on $L^2(G)$ is smoothing, since 
$\langle R_\theta z R_{-\theta}; 0\leq\theta\leq 1\rangle =G$.
In fact, geometrically if we identify $G$ with the unit sphere in $\Bbb C^2$ through the map $g=\pmatrix v&-\bar w\\ w&\bar v\endpmatrix \to\pmatrix v\\
w\endpmatrix$, $S_z$ is obtained by a circular average with circle of radius $s$ centered at $\pmatrix rv\\ rw\endpmatrix$ in the plane spanned by $\pmatrix -w\\
v\endpmatrix$ and $\pmatrix iw\\ iv\endpmatrix$.

In view of (10.41), we obtain that $(10.42)<\delta_2^c$ and hence $\Vert Tf\Vert_2 <\delta_1^{\frac 18 cc'}$.
This proves (10.40).
\bigskip

Next we treat the general case, using the result for $d=2$.

\noindent
{\bf Proof of Lemma 10.35 for $G=SU(d), d>2$.}

Take a subgroup $H$ of $G$, $H\simeq SU(2)$, considering for instance the embedding $SU(2)\to SU(d):h\mapsto
\sum^2_{i, j=1} h_{ij}\ e_i\otimes e_j$.

Write
$$
\align
\Vert F_1 *F_2\Vert_1&= \int\limits_G \Big|\int\limits_G F_1 (g) F_2 (g^{-1} x) dg\Big|dx\\
&\int\limits_G\int\limits_G\Big[\int\limits_H\Big|\int\limits_H F_1(hg)  F_2 (g^{-1} h^{-1} yx) dh\Big|dy\Big]dxdg.\tag 10.44
\endalign
$$
Fixing $x, g\in G$, introduce the following functions $\vp_1, \vp_2$ on $H$
$$
\aligned
\vp_1(y)&=F_1(yg)\\
\vp_2(y)&= F_2 (g^{-1} yx)
\endaligned
\tag 10.45
$$
for which the expression [ \ ] in (10.44) becomes
$$
\int\limits_H \Big|\int\limits_H\vp_1(h) \vp_2(h^{-1} y) dh\Big| dy =\Vert\vp_1* \vp_2\Vert_{L^1(H)}.
$$
In order to reach the conclusion by applying Lemma 10.35 on $H$, it will suffice to bound $\Vert\vp_1* P_{\delta_2}\Vert_{L^2(H)}< \delta_1^{c_1}$, for some
$\delta_2>\delta_1$, $\log \frac 1{\delta_2} \sim \log\frac 1{\delta_1}$, in the mean over $x, g\in G$ (cf. (10.44)).
Thus what is needed is an estimate of the form
$$
\int\limits_G\int\limits_G\int_H\Big|\operatornamewithlimits\nint\limits_{B_H(1, \delta_2)} F_2 (g^{-1} y_1 yx) dy_1\Big|^2 dydxdg<\delta_1^c,
$$
hence
$$
\int\limits_G\int\limits_G \Big|\operatornamewithlimits\nint\limits_{B_H(1, \delta_2)} F_2 (g_1 yg_2) dy\Big|^2 dg_1 dg_2 < \delta_1^c\tag 10.46
$$
Here $B_H(1, \delta_1) =\{y\in H; \Vert 1-y\Vert < \delta_1\}$ and $\displaystyle \nint$ denotes the average.

Rewrite (10.46) as
$$
\iint\limits_{GG}  \qquad 
\ \ {\not\!\!\!\!\!\!\!\!\!\!\int\limits_{B_H(1, 2\delta_2)}} F_2(g_2) \overline{F_2(g_1 y g_1^{-1} g_2)} dy\, dg_1 dg_2.\tag 10.47
$$
Fix $y\in B_H(1, 2\delta_2)$, $\Vert 1-y\Vert> \delta_2^2$ \big(the contribution of $\Vert 1-y\Vert\leq \delta_2^2$ in
(10.47) is at most $O(\delta_2^3)$\big).
We obtain$$
\iint F_2(g) \overline{F_2(xg)} \, \eta (dx)dg\tag 10.48
$$
where $\eta$ is the image measure $\Phi[\lambda_G]$ under the map
$$
\Phi =\Phi_y: G\rightarrow G: g\mapsto gyg^{-1}.\tag 10.49
$$
Next, taking $\eta_*$ the image of $\eta$ under $x\to x^{-1}$
$$
|(10.48)|^2 \leq \int \Big|\int F_2(xg) \eta(dx) \Big|^2 dg\leq \int\Big|\int F_2(xg)(\eta*\eta_*)(dx)\Big|^2 dg,
$$
and, similarly, for $r\in\Bbb Z_+$
$$
|(10.48)|^{2r}\leq\int\Big|\int F_2 (xg) (\eta*\eta_*)^{(r)} (dx)\Big|^2dg \tag 10.50
$$
where $\eta^{(r)} =\underbrace{\eta*\cdots*\eta}_{r-\text{fold}}$.

Returning to $\Phi$, diagonalize $y=\pmatrix e^{i\theta}& 0\\ 0& e^{-i\theta}\endpmatrix, \delta_2^2 <|\theta|< 2\delta_2$.
Hence
$$
\Phi (y) =1+i\theta gvg^{-1}+O(\theta^2)\tag 10.51
$$
where $v=e_1 \otimes e_1 -e_2 \otimes e_2 \in su(d)$ and
$$
\Phi (g_1)^{-1} \Phi(g_2) \phi(g_3)^{-1}\cdots \Phi (g_{2r})= 1+i\theta (-g_1v g_1^{-1} +g_2 v{g_2^{-1}} \cdots + g_{2r} vg_{2r }^{-1})+O(\theta^2).
$$
For $r =r(d)$ large enough, the map
$$
(g_1, g_2, \ldots, g_{2r})\mapsto -g_1 vg_1^{-1} +g_2 vg_2^{-1}-\cdots + g_{2r} v g_{2r}^{-1}
$$ 
gives a smooth density on $su(d)$.
Hence $(\eta*\eta_*)^{(r)}$ is a smooth density on $B(1, C|\theta|)\subset G$ \big(with derivative estimates in terms of $\frac 1{|\theta|} <\frac
1{\delta_2^2}$\big).

Recall that $F_2$ satisfies $\Vert F_2*P_{\delta_1}\Vert_2< \delta_1^c$.
Taking for $\delta_2$ an appropriate power of $\delta_1$, we get from the preceding
$$
\align
(10.50) = & \Vert(\eta*\eta_*)^{(r)}* F_2\Vert_2<\\
&\Vert\big((\eta*\eta_*)^{(r)}*P_{\delta_1}\big) * F_2\Vert_2 +O\big(\delta_1\delta_2^{-2}\big)<\\
&\Vert (\eta*\eta_*)^{(r)}\Vert_1 \Vert P_{\delta_1} * F_2\Vert_2+ O\big(\delta_1\delta_2^{-2}\big)< \delta_1^c+O\big(\delta_1\delta_2^{-2}\big)<\delta_1^{c'}.
\endalign
$$

This proves (10.46) and Lemma 10.35. \qed

\bigskip

\centerline
{\bf \S11. Expansion in $SU(d)$\quad (2)}

It remains to show that there is no approximate group $H$ satisfying (10.11) - (10.15).

Since $H$ is a union of $\delta$-balls, (10.14) is equivalent to
$$
\nu^{(\ell)} (aH)>\delta^{0+}.\tag 11.1
$$
Recall that $\ell\sim\log \frac 1\delta$.
Writing for $k<\ell$
$$
\nu^{(\ell)} (aH)= \sum_x \nu^{(\ell-k)}(x) \nu^{(k)} (xaH)
$$
it follows from (11.1) that for some $x\in G$
$$
\nu^{(k)} (xaH)> \delta^{0+}
$$
and hence
$$
\nu^{(2k)} (H.H)>\delta^{0+}.\tag 11.2
$$
In particular, recalling Kesten's bound [K] for random walks on $F_2$ 
$$
|H.H \cap W_k|>  \sqrt 2^k \ \text { if } \ 2^k> \Big(\frac 1\delta\Big)^{0+}.\tag 11.3
$$
Denote $V\subset \Mat_{d\times d}(\Bbb C)$ the real vector space of anti-Hermitian matrices of zero trace \big(i.e. the Lie-algebra of $G=SU(d)$\big),
which is irreducible under the adjoint representation of $G$ and hence of its Zariski-dense subgroup $\Gamma$.

We make the following assumption on $\nu$.

\proclaim
{Assumption $(*)$}
There is $\omega> 0$ such that for a proper subspace $L$ of $V$ and  $k$ large enough, one has the estimate
$$
\nu^{(k)} [g\in G; g^{-1} Lg=L]< e^{-\omega k}.\tag 11.4
$$
\endproclaim

Let $a\in V\backslash \{0\}$ and $2^k>\big(\frac 1\delta\big)^{0+}$.
Consider the increasing subspaces $L_s$ of $V$ defined by
$$
\align
L_1&= \span [a]\\
L_{s+1}&= \span [g^{-1} L_s g; g\in H.H \cap W_k].
\endalign
$$
Taking $s$ such that $L_{s+1} =L_s=L$, we have that $g^{-1} Lg=L$ for $g\in H.H\cap W_k$ and (11.2), (11.4) imply that $L=V$.
Since $(H.H\cap W_k)^{(s)} \subset H^{2s}\cap W_{sk}$, we proved

Assume $(*)$.
Let $a\in V\backslash\{0\}$ and $L$ a proper subspace of $V$.
Then, for $2^k>\big(\frac 1\delta\big)^{0+}$, there is $g\in H'\cap W_k$ such that $g^{-1} ag\not\in L$.

Equivalently, if $a, b\in V\backslash \{0\}$, there is $g\in H'\cap W_k$ such that
$$
\Tr g^{-1} agb^*\not= 0.\tag 11.5
$$

Recalling that $g_1, g_2$ are algebraic, (11.5) implies the following quantitative statement as a consequence of the effective Nullstellensatz (see Theorem
5.1 in [BY] and the comment on its generalization to polynomials with coefficients in the ring of integers in a fixed number field $K$, $[K:\Bbb Q]<\infty)$.

\proclaim
{Lemma 11.6} Assume $(*)$.
If $2^k>\big(\frac 1\delta\big)^{0+}$ and $a, b\in V\backslash \{0\}$, there is $g\in H'\cap W_k$ such that
$$
|\Tr g^{-1} agb^*|\geq C^{-k}\Vert a\Vert \ \Vert b\Vert\tag 11.7
$$
where $C$ is some constant depending on the generators $g_1, g_2$ of $\Gamma$.
\endproclaim

Next, apply Corollary 10 to $H.H\cap W_k, 2^k> \big(\frac 1\delta\big)^{0+}$.
We obtain
$$
i\xi \in V, \xi_{jj} =0 (1\leq j\leq d), \Vert\xi\Vert =1\tag 11.8
$$
such that for $\delta_0>\delta_1>0$, $\log \frac 1{\delta_0} \sim\log \frac 1{\delta_1}\sim k$
$$
\dist \Big(1+t\sum \xi_{ij}\eta_{ij}(e_i\otimes e_j), H'\Big)< \delta_1^{1+\delta}\Vert\eta\Vert\tag 11.9
$$
whenever $1+\eta\in H', \Vert\eta\Vert <\delta_1$ and $t\in [0, \delta_0]$.

Note that $\eta+\eta^*=0(\Vert\eta\Vert^2)$ and $\Tr \eta=0(\Vert\eta\Vert^2)$ and hence there is an element
$a\in V$ such that $\Vert\eta-a\Vert\lesssim \delta_1\Vert\eta\Vert$.
Thus from (11.9)
$$
\dist \Big(1+t \sum \xi_{ij} a_{ij} (e_i\otimes e_j), H')< \delta_0^{1+\gamma}\Vert\eta\Vert \ \text { for } \ t\in [0, \delta_0].\tag 11.10
$$
We may further replace in (11.10) $a$ by any conjugate $g^{-1}ag$ for $g\in H'$.
Take \hfill\break
$k_0\sim \log \frac 1{\delta_0}$ and small enough to ensure that
$$
C^{-k_0} >\delta_0^{\gamma/2}\tag 11.11
$$
where $C$ is the constant from Lemma 11.6.

Applying Lemma 11.6 with $a$ as above and $b= i\bar\xi \in V$ (where $\bar\xi_{ij}= \xi_{ji}$) gives some $g\in H'\cap  W_{k_0}$ such that
$$
\max_{i, j} |(g^{-1} ag)_{ij} \xi_{ij}|\gtrsim \Big|\sum_{i, j} (g^{-1} ag)_{ij} \xi_{ij}\Big| > C^{-k_0}\Vert\eta\Vert.\tag 11.12
$$

Let $\zeta\in V, \Vert\zeta\Vert=1$, be defined by normalization of  $\big((g^{-1} ag)_{ij} \xi_{ij}\big)_{1\leq i, j\leq d}$.

Clearly, from (11.10) and the preceding
$$
\dist (1+t\zeta, H')< \delta_0^{1+\gamma}\Vert \eta\Vert \ \text { for } \ t\in [0, \delta_0C^{-k_0} \Vert\eta\Vert]\overset {(11.11)}\to \supset
[0, \delta_0^{1+\frac \gamma 2} \Vert\eta\Vert].\tag 11.13
$$
Again from Lemma 11.6, there are elements $g_1, \ldots, g_{d_1} \in H'\cap W_{k_0} (d_1=d^2-1)$ such that
$$
|\det (g_s^{-1} \zeta g_s; 1\leq s\leq d_1)| > C_1^{-k_0}.\tag 11.14
$$
Since in (11.13) we may replace $\zeta$ by conjugates $g^{-1}\zeta g$ with $g\in H'$, it
easily follows from (11.14) that
$$
\dist (1+tV, H')< \delta_0^{1+\gamma}\Vert\eta\Vert
$$
for
$$
t\in [0, \delta_0 C_1^{-k}\Vert\eta\Vert]\supset [0, \delta_0^{1+\frac \gamma 2}\Vert\eta\Vert].\tag 11.15
$$
Hence, we proved (redefining $\delta_0$ and $\gamma$).

\proclaim
{Lemma 11.16} Let $1+\eta \in H', \Vert\eta\Vert<\delta_1$. Then
$$
\dist (1+tV, H')<\delta_0^{1+\gamma}\Vert\eta\Vert \ \text { for } \ t\in [0, \delta_0\Vert\eta\Vert].\tag 11.17
$$
\endproclaim

Fix a small constant $\ve$ and let $k=[\ve\log\frac 1\delta]$.
Thus $\delta_0 =\delta^{\ve_0}, \delta_1=\delta^{\ve_1}$ with $\ve_0\sim\ve_1\sim \ve$ in Lemma 11.16.

The final step consists in using Lemma 11.16 to derive a contradiction on (10.16), i.e.
$$
|H'|<\delta^{\gamma-}.\tag 11.18
$$
Fix an element $g_1 =1+\eta_1\in H.H\cap W_{k_1} \cap B(1, \delta_1), g_1\not= 1$, which by (11.3) is possible for $k_1<\frac 1{d_1} \log\frac
1{\delta_1}$.
It follows that
$$
\delta_1^C<\Vert\eta_1\Vert = t_1< \delta_1.\tag 11.19
$$
From (11.17)
$$
\dist (1+tV, H')<\delta_0^{1+\gamma} t_1 \ \text { for } \ t\in [0, \delta_0t_1].\tag 11.20
$$
Hence, for any $t_2 \in [\delta_0^{1+\frac\gamma 2} t_1, \delta_0t_1]$, $H'$ contains some element $1+\eta$ with $\Vert\eta\Vert\approx t_2$.
Applying again Lemma 11.16 shows that
$$
\dist (1+tV, H')< \delta_0^{1+\gamma} t_2 \ \text { for } \ t\in [0, \delta_0 t_2]\tag 11.21
$$
and therefore $H'$ contains elements $1+\eta$ with $\Vert\eta\Vert\approx t_3$ for any $t_3 \in [\delta_0^{1+\frac \gamma 2}t_2, \delta_0 t_2]$
and therefore any $t_3 \in [\delta_0^{2(1+\frac\gamma 2)} t_1, \delta_0^2 t_1]$.

After $r$ steps, we see that $H'$ contains elements $1+\eta$ with $\Vert\eta\Vert\approx t$ for any $t\in [\delta_0^{r(1+\frac \gamma 2)} t_1, \delta_0^r
t_1]$.
Taking $r=[\frac 2\gamma]+1$, it follows that $H'$ contains elements $1+\eta$ where $\Vert\eta\Vert \approx t$, for any $t\in [\delta, \delta_0^r t_1]$.

Another application of Lemma 11.16 implies
that
$$
\dist (1+tV, H')< \delta_0^\gamma t \ \text { for } \ t\in [0,\delta_2]\tag11.22
$$
where
$$
\delta_2=\delta_0^{r+1} t_1.\tag 11.23
$$
We claim that
$$
G\cap B\big(1, \delta_2\big) \subset H'.\tag 11.24
$$
Then
$$
|H'|> \delta_2^{d_1} > \delta^{(\ve_0(r(\gamma)+1)+C\ve_1)d_1}
$$
contradicting (11.18) for $\ve$ small enough.

\noindent
{\bf Proof of (11.24)}  Take $g_0\in G\cap B(1, \delta_2)$.
Then $g_0=1+\eta_0$ and there is $a_0\in V$ such that $\Vert a_0-\eta_0\Vert\lesssim \delta_2^2$.
By (11.22), $\dist (1+a_0, H')<\delta_0^\gamma\delta_2$ and we take $h_1\in H'$ such that
$$
\Vert g_0h_1^{-1}-1\Vert =\Vert g_0-h_1\Vert\lesssim \delta_2^2 +\delta_0^\gamma\delta_2 \lesssim \delta_0^\gamma\delta_2.\tag 11.25
$$
Next, write $g_1 =g_0 h_1^{-1} = 1+\eta_1$ and take $a_1 \in V$ with $\Vert a_1-\eta_1\Vert\lesssim \Vert a_1\Vert^2$.
By (11.22), $\dist (1+a_1, H')< \delta_0^\gamma\Vert a_1\Vert$ and we obtain $h_2 \in H'$ such that
$$
\Vert g_0h_1^{-1} h_2^{-1}-1\Vert =\Vert g_1-h_2\Vert\lesssim \Vert a_1\Vert^2 +\delta_0^\gamma\Vert a_1\Vert\lesssim \delta_0^{2\gamma}
\delta_2.\tag 11.26
$$
Since $H'$ is a union of $\delta$-balls, a few iterations give the desired conclusion.

This completes the proof of the spectral gap, conditional  on the assumption $(*)$.

\bigskip

\centerline
{\bf \S12. Proof of Assumption $(*)$ for $d=3$}

Assume $H\subset \Gamma$ and $L$ a nontrivial subspace of $V=su(d)$ satisfying

\roster
\item "{(12.1)}" $\nu^{(k)}(H)> e^{-\ve k}$

\item "{(12.2)}" $g^{-1} Lg =L$ for $g\in H$.
\endroster

\noindent
where in (12.1) we assume $\ve$ a sufficiently small constant (depending on $\Gamma$ and $\nu$) and $k$ large.
Our purpose is to get a contradiction for $d=3$.
This will illustrate the method in the simplest case.
The argument in the general case is given in \S14.
Essential use is made of the theory of random matrix products as developed by Furstenberg and Guivarch.
A treatment of this theory in the setting of general local fields appears in [A].

Recall that $\Gamma \subset SU(k)$, where $k$ is the algebraic closure of $\Bbb Q$.
We will consider $V$ and $L$ as vector spaces over $k$, hence $V=\{g\in\Mat_{d\times d}(k), \Tr g=0\}$.

\noindent
{(i)} Exploiting the theory of random matrix products requires proximal elements.
Following the approach of Tits [Tits]  (see also [G]) proximal elements in a suitable setting may be produced by passing to an appropriate local field.

Fix an element $g_0\in\Gamma$ with eigenvalues $(\lambda_j)_{1\leq j\leq d}$ such that  not all quotients $\frac{\lambda_j}{\lambda_{j'}}$ are roots of
unity.
If $\frac{\lambda_j}{\lambda_{j'}}$ is not root of unity, there is a local field $k\subset K_v$ such that $v(\lambda_j/\lambda_{j'})\not= 1$.
Hence we have
$$
|\{v(\lambda_j); 1\leq j\leq d\}|\geq 2.
$$
For $d=3$, either
$$
v(\lambda_1) > v(\lambda_2)> v(\lambda_3)\tag 12.3
$$
or
$$
v(\lambda_1)=v(\lambda_2)> v(\lambda_3)\tag 12.4
$$
(if $v(\lambda_1)>v(\lambda_2) =v(\lambda_3)$, replace $g_0$ by $g_0^{-1}$).

Denote $\rho$ the adjoint representation on $V$.

If (12.3) holds, the representation $\rho|\Gamma$ on $V\otimes K_v$ has $\rho_{g_0}$ as proximal element and since it is totally irreducible
(by the Zariski density assumption), random matrix product theory implies that (12.1), (12.2) are not compatible for $\ve>0$ small enough.

Thus we may assume in the sequel that the situation (12.3) may not be realized for any $g\in\Gamma$.

\noindent
(ii) Consider the representation of $\Gamma$ on $\overset 2\to \Lambda (V\otimes K_v)$, which we also denote $\rho$.

If $e_1, e_2, e_3$ diagonalizes $g_0, g_0 e_i=\lambda_i e_i (1\leq i\leq 3)$, the eigenvector
$$
\xi = (e_3 \otimes e_1) \wedge (e_3\otimes e_2)\tag 12.5
$$
of $\rho_{g_0}$ has dominant eigenvalue $\frac{\lambda_1^2}{\lambda_3^2}$, by (12.4).

Denote by 
$$
\Cal S=\span_k [\rho_g(\xi); g\in\Gamma]\tag 12.6
$$
the subspace of $\overset 2\to \Lambda V$.
The restriction of $\rho$ to $\Cal S$ is totally irreducible.
Otherwise, there would be a proper subspace $\Cal S_1$ of $\Cal S$ which is invariant under a finite index subgroup $\Gamma_1$ or $\Gamma$.
Hence $\Cal S_1\otimes \Bbb C$ would be invariant for $\rho_g, g\in\bar\Gamma_1= $ Zariski closure of $\Gamma_1$.
Since $\bar\Gamma = SL_d(\Bbb C)$ and $\bar \Gamma_1$ is a finite index subgroup of $\bar \Gamma, \bar\Gamma_1 =SL_d(\Bbb C)$.
In particular $\Cal S_1$ is $\Gamma$-invariant, hence $\Cal S_1=\Cal S$. Also
$$
\Cal S\otimes \Bbb C=\span_{\Bbb C} [\rho_g(\xi); g\in SL_d(\Bbb C)].\tag 12.7
$$
Since $\rho$ restricted to $\Cal S\otimes K_v$  has a proximal element, it follows again from random matrix product theory and (12.1),  that
$$
\Cal S=\span[\rho_g(\eta); g\in H] \ \text { for any } \ \eta\in \Cal S\backslash\{0\}.\tag 12.8
$$
We used here that the probabilistic estimates depend on $\nu$ but not on the vector $\eta$.

\noindent
{(iii)} The space $\overset 2\to\Lambda V$ decomposes as
$$
\overset 2\to \Lambda V =\overset 2\to \Lambda L \ \underline\oplus \ \overset 2\to \Lambda L^\bot\, \underline\oplus\,
(L\wedge L^\bot) =\frak S_1  \ \underline\oplus \ \frak S_2 \ \underline\oplus  \ \frak W\tag 12.9
$$
and by (12.2), each of the components is invariant under $\rho_g$ for $g\in H$.
Take an element $g_1\in H$ such that not all quotients of its eigenvalue are roots of unity
(this is certainly possible, since $\log |H\cap W_k|\sim k$).
Arguing as in (i) and since case (12.3) was ruled out, we are in the situation (12.4) (in some local field $K_w)$.
Consider the representation on $\overset 2\to\Lambda(V\otimes K_w)$ as in (ii).
Note that if $\frak X$ is a subspace of $\overset 2\to\Lambda V$ invariant under $\rho_{g_1}$, then its eigenvector 
$(e_3'\otimes e_1')\wedge (e_3' \otimes e_2')$ with top exponent will either be orthogonal on $\frak X$ or belong to $\frak X\otimes K_w$, hence to $\frak X$.

Therefore, considering the decomposition (12.9), it follows that $(e_3'\otimes e_1')\wedge (e_3'\otimes e_2')=\eta$ belongs to one of the spaces $\frak
S_1, \frak S_2$ or $\frak W$.
Also, by (12.7), $\Cal S$ contains any element of the form $(x\otimes y)\wedge (x\otimes z)$ with $x, y, z \in k^3$ and $\langle x, y\rangle =0=\langle x,
z\rangle$.
In particular $\eta\in\Cal S$ and it follows from (12.8) that $\Cal S$ is contained in one of the spaces $\frak S_1, \frak S_2$ or $\frak W$.

There are now three cases to conisder.

\noindent
{\bf Case I.} $\Cal S\subset \overset 2\to\Lambda L$.

Since $(x\otimes y)\wedge(x\otimes z)\in\overset 2\to \Lambda L$ for all $x, y, z\in k^3, \langle x, y\rangle=0 =\langle x, z\rangle$, it follows that
$x\otimes y \in L$ whenever $x, y\in k^3, \langle x, y\rangle =0$.
Therefore $L=V$, a contradiction.

\noindent
{\bf Case II.}
$\Cal S\subset\overset 2\to \Lambda L^\bot$.  The same argument as in case I applies.

\noindent
{\bf Case III.} $\Cal S\subset L\wedge L^\bot$.

Note that if $a, b\in V$ and $a\wedge b\in L\wedge L^\bot$, then $L$ and $L^\bot$ both contain a nontrivial linear combination of $a$ and $b$.

Assume $\dim L\leq \dim L^\bot$, hence $\dim L\leq 4$.

Considering the sectors $(x\otimes y)\wedge (x\otimes z)\in L\wedge L^\bot$, it follows from the preceding that for each $x\in k^3$, there is some
$x'\not=0, \langle x, x'\rangle =0$ such that
$$
x\otimes x' \in L.\tag 12.10
$$
Our next aim is to show that (12.10) forces $\dim L\geq 5$, hence again a contradiction.

Consider the real algebraic variety
$$
\Omega =\{(x, y) \in\Bbb C^3 \times \Bbb C^3|\langle x, y\rangle =0 \ \text { and } \ x\otimes y\in L\otimes \Bbb C\}\tag 12.11
$$
(an intersection of quadrics).

By (12.10), we may introduce a real-analytic function $\vp:O\to \Bbb C^3\backslash\{0\}$, $O\subset \Bbb C^3$ some open set, such that for all
$x\in O$

\roster
\item"{(12.12)}" $\langle x, \vp(x)\rangle = 0$.
\item "{(12.13)}"  $x\otimes \vp(x) \in L\otimes \Bbb C$.
\endroster

We distinguish two further cases. 

\noindent
{(a)}. $\vp$ has 2-dim range (over $\Bbb C$).

If $\Im \vp \subset [e_1, e_2]$, then necessarily, for $x=x_1 e_1+x_2e_2+x_3e_3$,  $\vp(x)$ is parallel to $ -\bar x_2 e_1+\bar x_1 e_2$ by (12.12), implying
$$
(x_1 e_1+x_2e_2+x_3e_3)\otimes (-\bar x_2 e_1+\bar x_1 e_2)\in L\otimes \Bbb C\tag 12.14
$$
for all $x\in O$ and hence for all $x\in\Bbb C^3$.

Since the functions $x_1^2, x_2^2, x_1x_2, x_1x_3, x_2x_3$ are linearly independent,
$$
e_1\otimes e_2, e_2\otimes e_1, e_1 \otimes e_1-  e_2\otimes e_2,  e_3\otimes e_2, e_3\otimes e_1 \in L
$$
and $\dim L\geq 5$.

(b) $\vp$ has 3-dim range.

We may clearly find elements $x_1, x_2, x_3 , x_4, x_5 \in O$
such that for each triplet $\{i, j, k\} \subset \{1, 2, 3, 4, 5\}$ of distinct integers, each of the systems $\{x_i, x_j, x_k\}$ and
$\{\vp(x_i), \vp(x_j), \vp(x_k)\}$ consist of linearly independent vectors.

WE claim that $\{x_i\otimes \vp (x_i); i=1, \ldots, 5\}$ are linearly independent, which can be  seen as follows.
Fix an index $i=1$.
From our assumptions, there is $T\in\Cal L(\Bbb C^3, \Bbb C^3)$ such that $Tx_2= Tx_3=0$ and $\langle Tx_1, \vp (x_1)\rangle \not= 0$,
$\langle Tx_1, \vp(x_4)\rangle =0= \langle Tx_1, \vp (x_5)\rangle$.
Hence 
$$
\langle Tx_2, \vp(x_2)\rangle =0= \langle Tx_3, \vp(x_3)\rangle $$
and writing $$ x_4=a_4x_1+b_4x_2+c_4x_3,$$
$$x_5 =a_5x_1+ b_5x_2+c_5 x_3,$$  we get $$ \langle Tx_4, \vp (x_4)\rangle =a_4 \langle T x_1, \vp(x_4)\rangle = 0= \langle Tx_5, \vp (x_5)\rangle.$$

This completes the proof of the main theorem for $SU(3)$.

\bigskip

\centerline
{\bf \S13. Lemmas on linear independence}
\proclaim
{Lemma 13.1} Let $1\leq k\leq d$ and $\vp_1, \ldots, \vp_k$ be continuous complex, linearly independent functions on $\Bbb
C^d$. Then
$$
\dim [x_i\vp_j(x); 1\leq i\leq d \ \text { and } \ 1\leq j\leq k]\geq k\Big(d-\frac{k-1}{2}\Big).\tag 13.2
$$
where $x_i\vp_j(x)$ are viewed as functions on $\Bbb C^d$.
\endproclaim

\noindent
{\bf Proof.}

The proof is by induction on $d$.

Denote $e_1, \ldots, e_d$ the unit vector basis of $\Bbb C^d$.

Take $k_1\leq k$.
Since $\vp_1, \ldots, \vp_{k_1}$ are linearly independent functions, there are $\xi_1, \ldots, \xi_{k_1} \in\Bbb C^d$ s.t.
$$
\det [\vp_j(\xi_{j'})]_{1\leq j, j'\leq k_1} \not= 0.\tag 13.3
$$
By a linear transformation, we may assume that $\xi_1, \ldots, \xi_{k_1} \in [e_1, \ldots, e_{k_1}]$ and hence
$\vp_j|_{[e_1, \ldots,
e_{k_1}]}(1\leq j\leq k_1)$ are linearly independent.

We distinguish 2 cases.

\noindent
{\bf Case 1:} $k=d$

Taking $k_1=d-1$, we can assume that $\vp_j|_{[e_1, \ldots, e_{d-1}]} (1\leq j\leq d-1)$ are linearly independent.
Denote $x' =x_1' e_1+\cdots+ x_{d-1}' e_{d-1}$.

From the induction hypothesis, there is a subset $\Omega\subset \{ (i, j); 1\leq i, j\leq d-1\}$ such that $|\Omega|\geq
\frac {d(d-1)}2$
and $\big(x_i' \vp_j(x')\big)_{(i, j)\in \Omega}$ are linearly independent functions in $x'$.  We claim that
$$
\big(x_i \vp_j(x)\big)_{(i, j)\in\Omega}\cup \big(x_d\vp_j(x)\big)_{1\leq j\leq d}
$$
are linearly independent, which will imply that $\dim [x_i\vp_j (x)]\geq \frac {d(d-1)}{2}+d$.

Suppose the claim fails.
Then there is a nontrivial linear combination
$$
\sum_{(i, j)\in\Omega} a_{ij} x_i\vp_j(x)+\sum^d_{j=1} a_{d_j} \ x_d\vp_j(x)=0.
$$
Setting $x_d=0$, we get
$$
\sum_{(i, j)\in\Omega} a_{ij}x_i' \vp_j(x')=0  \ \text { hence } \ a_{ij} =0\ \text { for } \ (i, j)\in\Omega.
$$
Therefore $x_d\sum^d_{j=1} a_{dj} \vp_j(x) =0$ and $\sum^d_{j=1} a_{dj} \vp_j(x)=0$, since the $\vp_j$ are continuous.
Hence, also $a_{dj}=0$  \ (contradiction).

\noindent
{\bf Case 2: $k< d$}

Take $k_1 =k$ and argue as above to obtain
$$
\dim [x_i\vp_j(x)]\geq k\Big( d-1-\frac {k-1}2\Big) +k = k\Big(d-\frac{k-1}{2}\Big).
$$
This proves the Lemma.\qed

\noindent
{\bf Remark.} The assumption that the $\vp_j$ are continuous in Lemma 13.1 can not be dropped.
Take, for instance,  a basis $e_1, \ldots, e_d$ and define
$$
\left\{
\aligned &{\vp_j (e_j)=1}\\
&{\vp_j(x) =0 \text { if } x\not= e_j.}
\endaligned
\right.
$$
Since $x_i\vp_j=\delta_{ij}\vp_j$, $\dim [x_i\vp_j; 1\leq i, j \leq d]=d$.\qed

\proclaim
{Lemma 13.4}
Let $1\leq k\leq d-1$ and $\vp: O\to G_{d, k}, O\subset\Bbb C^d$ some open set, a continuous map, satisfying
$\vp(x) \subset [x]^\bot$.
Then
$$
\dim [x\otimes\vp(x); x\in O]\geq (k+1) d-1.\tag 13.5
$$
\endproclaim

\noindent
{\bf Proof.}
One may clearly choose a subspace $E$ of $\Bbb C^d$ such that $dim E= d-k+1$ and

(13.6) $\dim [\vp(x) \cap E]=1$.

(13.7) $\Proj_{E^\bot}\vp(x)=E^\bot$.

\noindent
for $x\in O$.
Hence
$$
\dim [x\otimes\vp (x)]\geq \dim [x\otimes \big(\vp(x)\cap E)]+\dim [x\otimes P_{E^\bot} \vp(x)]=(13.8)+d(k-1).
$$
Introduce a continuous function $\psi: O\to \Bbb C^d\backslash\{0\}$ such that $\psi(x) \in\vp(x)\cap E$, hence
$\langle x, \psi(x)\rangle =0$.
Since clearly $\dim [\psi_1, \ldots, \psi_d]\geq 2$, application of Lemma 13.1 with $k=2$ gives
$$
(13.8)\geq \dim [x\otimes \psi(x)]\geq 2k-1
$$
This proves Lemma 13.3.

\bigskip
\centerline
{\bf \S14. Assumption $(*)$ \quad (General Case)}

Following preceding analysis for $d=3$, we may introduce the set

$\Cal D =\{(d_+, d_-) \in \{1, \ldots, d-1\}^2$; there is $g\in\Gamma$ and a local field $K_v$ such that the exponents of
$g$ may be ordered as
$$
v(\lambda_1) =\cdots = v(\lambda_{d_+}) > v(\lambda_{d_++1})\geq \cdots > v(\lambda_{d-d_ - +1 }) =\cdots =
v(\lambda_d\}
\tag 14.1
$$
where $d_+, d_-< d$, as we assume $v(\lambda_i)(1\leq i\leq d)$ not all equal and $d_++d_-\leq d$.

Fixing a configuration $(d_+, d_-)\in \Cal D$, we obtain a proximal representation by considering the extension of the
adjoint
representation to the exterior power
$$
\overset D\to \Lambda(V\otimes K_v)\tag 14.2
$$
where $D=d_+. d_-$.
The proximal vector is given by
$$
\xi =\operatornamewithlimits\Lambda\limits_{\Sb 1\leq i\leq d_+\\ d-d_-< j\leq d\endSb} (e_i\otimes e_j)\tag 14.3
$$
in a suitable OB $\{e_1, \ldots, e_d\}$; the eigenvalue is $(\frac {\lambda_1}{\lambda_d})^D$.

Denote
$$
\Cal S=\span_k [\rho_g(\xi); g\in \Gamma]\tag 14.4
$$
the subspace of $\overset D\to \Lambda V$.
Again from Zariski density of $\Gamma$ in $SL_d(\Bbb C)$,
$$
\Cal S\otimes \Bbb C =\span_{\Bbb C} [\rho_g(\xi); g\in GL_d (\Bbb C)]\tag 14.5
$$
and $\Gamma$ acts strongly  irreducibly on  $\Cal S$.

From random matrix product theory, also
$$
\Cal S=\span[\rho_g(\eta); g\in H] \ \text { for any } \ \eta\in\Cal S\backslash \{0\}\tag 14.6
$$
provided $H$ satisfies
$$
\nu^{(k)} (H) > e^{-\ve k}\tag 14.7
$$
with $k$ large enough and $\ve$ small enough.

\bigskip
{\bf (2)}. Assume given a nontrivial subspace $L$ of $V$ satisfying
$$
\rho_g(L) =L \ \text { for } \ g\in H.\tag 14.8
$$
The space $\overset D\to\Lambda V$ decomposes as the direct sum
$$
\overset {D_0} \to\Lambda L \wedge \overset{D_1}\to \Lambda L^\bot\tag 14.9
$$
where $D_0+D_1 =D$.
Note also that since $L$ was obtained as complexification of a subspace of $su(d)$,  we have that $L=L^*$.

Taking some element $g\in H$ which  has the property that is eigenvalue quotients are not all roots of unity and considering an appropriate
valuation,
we obtain some type $(d_+, d_-)\in\Cal S$ with expanding vector $\eta\in\overset D\to \Lambda V$ of the type (14.3).
Note that we may always assume that $d_+\geq d_-$ since $g$ may be replaced by $g^{-1}$.

Since the components of the decomposition (14.9) are $g$-invariant, we conclude that
$$
\eta\in \overset {D_0}\to \Lambda L\wedge \overset {D_1}\to \Lambda L^\bot\tag 14.10
$$
for some $D_0, D_1, D=D_0+D_1$.
By (14.6)
$$
\Cal S\subset \overset {D_0}\to\Lambda L\wedge \overset {D_1}\to \Lambda L^\bot.\tag 14.11
$$
We also note that from (14.5), $\Cal S\otimes \Bbb C$ contains any element of the form
$$
\Lambda_{\Sb 1\leq i\leq d_+\\ 1\leq j\leq d_-\endSb} (x_i\otimes y_j)\tag 14.12,
$$
where $\{x_1, \ldots, x_{d_+}, y_1, \ldots, y_{d_-}\}$ are orthogonal vectors in $\Bbb C^d$ (for the Hermitian inner
product).

By (14.11), it follows that $L$ (respectively  $L^\bot$) will contain $D_0$ (respectively $D_1$) linearly independent elements from
$$
\span [x_i\otimes y_j; 1\leq i\leq d_+, 1\leq j\leq d_-].\tag 14.13
$$
Note that if $D_1=0$, then obviously $L$ contains (14.13) and hence any element $x\otimes y$ with $\langle x, y\rangle =0$.
Thus $L=V=\{ x\in\Mat_{d\times d} (\Bbb C); \Tr\, x=0\}$ would be trivial.

Hence we assume $D_0\geq 1, D_1\geq 1$.

It follows from (14.13) that, given orthogonal subspaces
$E_+, E_-$ of $\Bbb C^d$, $\dim E_+=d_+, dim E_-=d_-$,
$$
\dim \big(L\cap (E_+ \otimes E_-)\big)+\dim \big(L^\bot\cap (E_+\otimes E_-)\big) =\dim (E_+\otimes E_-).\tag 14.14
$$
Hence
$$
\dim \Proj_L(E_+\otimes E_-)+\dim\Proj _{L^\bot} (E_+\otimes E_-)=\dim (E_+\otimes E_-).\tag 14.15
$$
Denoting $F_0=\Proj_L(E_+\otimes E_-), F_1=\Proj_{L^\bot}(E_+\otimes E_-)$, clearly $E_+\otimes E_-\subset F_0+F_1$ and
(14.15) implies $E_+\otimes E_- = F_0+F_1$.
Therefore
$$
\Proj _L (E_+\otimes E_-)\subset E_+\otimes E_-, \Proj_{L^\bot}(E_+\otimes E_-)\subset E_+\otimes E_-.\tag 14.16
$$
Next, let $x, y\in\Bbb C^d\backslash \{0\}, \langle x, y\rangle = 0$.
From (14.16)
$$
\Proj _L (x\otimes y)\in \bigcap_{\Sb x\in E_+\\ y \in E_-\endSb} (E_+\otimes E_-)\equiv S_{x, y} \subset V.\tag 14.17
$$
Assume $d_+\geq 2, d_-\geq 2$.

We claim that $S_{x, y} =[x\otimes y]$.
For if $T\in S_{x, y}$, we have
$$
\text{Im\,} T\subset\bigcap E_-\tag 14.18,
$$
where $E_-$ ranges over all $d_-$-dimensional subspaces of $[x]^\bot$ such that $y\in E_-$.
Since $d_+ \geq 2, d_- < d-1$ and $(14.18) =[y]$.

Similarly, since $T^*\in \bigcap_{\Sb x\in E_+\\ y\in E_-\endSb} (E_-\otimes E_+)$, it follows that
$$
({\text{\rm Ker\,}} T)^\bot =\Im T^* \subset [x].
$$
Hence $T\in [x\otimes y]$, proving the claim.

Thus
$$
\Proj_L (x\otimes y)\in [x\otimes y] \ \text { and } \ \Proj_{L^\bot}(x\otimes y)\in [x\otimes y]
$$
implying that
$$
x\otimes y \in L \ \text { if } \ \Proj_L (x\otimes y)\not= 0\tag 14.19
$$
and similarly for $L^\bot$.

Fixing orthogonal vectors $e, e' \in \Bbb C^d\backslash\{0\}$, either $\Proj_L(e\otimes e')\not= 0$ or $\Proj_{L^\bot} (e\otimes
e')\not= 0$.
If $\Proj_L(e\otimes e') \not=0$, clearly $\Proj_L(x\otimes y) \not= 0$ for $x\in U, y\in U'$, $\langle x, y\rangle =0$, with $U$ (resp.
$U'$) some neighborhood of $e$ (resp $e'$).
From (14.19)
$$
(U\otimes U')\cap V\subset L\tag 14.20
$$
which is easily seen to imply that $V=L$ (contradiction).

It remains to consider the case $d_+=1$ (and similarly $d_-=1$).

Taking $E_+=[x], x\in \Bbb C^d\backslash \{0\}$, it follows from (14.14) that given any subspace $E_-$ of $[x]^\bot, \dim E_-=d_-$, 
there is a decomposition $E_-=W_0+W_1$ such that $x\otimes W_0\subset L, x\otimes W_1\subset L^\bot$.
Therefore clearly, for given $x\in\Bbb C^d\backslash \{0\}$
$$
\dim \big(L\cap (x\otimes[x]^\bot)\big)+\dim \big(L^\bot\cap (x\otimes [x]^\bot)\big) =d-1.\tag 14.21
$$
Specifying $k_0=\dim\big(L\cap (x\otimes [x]^\bot)\big)$, $k_1=\dim \big(L^\bot\cap (x\otimes[x]^\bot)\big)$ for $x$
restricted to some open subset $O\subset \Bbb C^d$, Lemma 13.4 in \S13 implies
$$
\dim L\geq (k_0+1) d-1 \ \text { and } \ \dim L^\bot\geq (k_1+1)d-1
$$
and hence, by (14.21)
$$
d^2-1\geq (d+1) d-2\quad 
$$
obtaining a contradiction and completing  the proof of the main theorem.

\enddocument

%% file: fig1.tex
\input pictex
\font\thinlinefont=cmr5
$$
\hbox{\beginpicture
\setcoordinatesystem units <.80000cm,0.80000cm>
\setshadesymbol ({\thinlinefont .})
\setlinear
%
%
\linethickness= 0.500pt
\setplotsymbol ({\thinlinefont .})
{\plot  2.540 20.796  9.525 25.400 / }%
%
%
\put{$\bullet$} [lB] at  2.500 20.700
%
%
\put{$\bullet$} [lB] at  3.728 21.549
%
%
\put{$\bullet$} [lB] at  5.499 22.701
%
%
\put{$\bullet$}  [lB] at  9.395 25.291
%
%
\put{$y_{j_1, \ldots, j_s, j_{s+1}}$} [lB] at 10.319 25.291
%
%
\put{$y_{j_1, \ldots, j_{s}, 1}$} [lB] at  4.604 21.590
%
%
\put{$y_{j_1,\ldots, j_s}$ } [lB] at  3.175 20.637
\linethickness=0pt
\putrectangle corners at  2.515 25.521 and 10.319 20.523
\endpicture}
$$

%% file: fig2.tex
\font\thinlinefont=cmr5
$$
\hbox{\beginpicture
\setcoordinatesystem units <.80000cm,.80000cm>
\setshadesymbol ({\thinlinefont .})
\setlinear
%
%
\linethickness= 0.500pt
\setplotsymbol ({\thinlinefont .})
{\putrule from  1.905 21.590 to 18.574 21.590
}%
%
%
\linethickness= 0.500pt
\setplotsymbol ({\thinlinefont .})
{\putrule from  1.905 19.050 to 18.574 19.050
}%
%
%
\linethickness= 0.500pt
\setplotsymbol ({\thinlinefont .})
{\putrule from  3.400 19.050 to  3.400 23.019 
} %
\linethickness= 0.500pt
\setplotsymbol({\thinlinefont .})
\
\linethickness= 0.500pt
\setplotsymbol ({\thinlinefont .})
\setdashes < 0.1270cm>
{\plot  3.135 19.051  9.525 24.130 / }
%
%
\linethickness= 0.500pt
\setplotsymbol ({\thinlinefont .})
{\plot  8.096 20.003  9.525 23.971 / }%
%
%
\linethickness= 0.500pt
\setplotsymbol ({\thinlinefont .})
{\plot 10.160 23.654 16.510 19.844 /
}%
%
%
\put{$\bullet$} [lB] at  2.820 18.924
%
%
\put{$\eta$} [lB] at  0.953 21.431
%
%
\put{$\bullet$ } [lB] at  2.820 21.480
\put{$0$} [lB] at 2.820  18.500 
%
\put{$\bullet$} [lB] at  7.960 19.903
%
%
\put{$10Kz_1$} [lB] at  9.850 23.872
%
%
\put{$\bullet$} [lB] at 16.440 19.675
%
%
\put{$\bullet$} [lB] at  6.120 21.480
%
%
\put{$\bullet$} [lB] at  9.575 24.109
%
%
\put{$z_3$} [lB] at 16.898 19.785
%
%
\put{$z_2$ } [lB] at  8.404 19.744
%
%
\put{$z_1$} [lB] at  6.600 20.999
\linethickness=0pt
\putrectangle corners at  0.953 24.162 and 18.599 18.574
\endpicture}
$$